\documentclass[11pt]{article}
\usepackage{latexsym,amssymb}
\usepackage[abbrev]{amsrefs}
\usepackage{color}
\usepackage{amsmath}

\newtheorem{Theorem}{\bf Theorem}[section]
\newtheorem{Lemma}{\bf Lemma}[section]
\newtheorem{Proposition}{\bf Proposition}[section]
\newtheorem{Corollary}{\bf Corollary}[section]
\newtheorem{Remark}{\bf Remark}[section]
\newtheorem{Example}{\bf Example}[section]
\newtheorem{Definition}{\bf Definition}[section]

\newenvironment{theorem}{\begin{Theorem}$\!\!\!$}{\end{Theorem}}
\newenvironment{lemma}{\begin{Lemma}$\!\!\!$}{\end{Lemma}}

\newenvironment{corollary}{\begin{Corollary}$\!\!\!$}{\end{Corollary}}
\newenvironment{remark}{\begin{Remark}$\!\!\!$}{\end{Remark}}

\numberwithin{equation}{section}

\begin{document}

\title{Refined asymptotic expansions of solutions\\ to fractional diffusion equations}
\author{Kazuhiro Ishige and Tatsuki Kawakami}
\date{}
\maketitle
\begin{abstract}
In this paper, as an improvement of the paper [K. Ishige, T. Kawakami and H. Michihisa, SIAM J. Math. Anal. 49 (2017) pp. 2167--2190], 
we obtain the higher order asymptotic expansions of the large time behavior of 
the solution to the Cauchy problem for inhomogeneous fractional diffusion equations 
and nonlinear fractional diffusion equations. 
\end{abstract}
\vspace{10pt}
\vspace{3pt}

\vspace{40pt}
\noindent Addresses:


\smallskip
\noindent 
K.~I.: \,Graduate School of Mathematical Sciences, The University of Tokyo,\\
\qquad\quad 3-8-1 Komaba, Meguro-ku, Tokyo 153-8914, Japan. \\
\noindent 
E-mail: {\tt ishige@ms.u-tokyo.ac.jp}\\

\smallskip
\noindent 
T. K.: Applied Mathematics and Informatics Course,\\
\qquad\quad Faculty of Advanced Science and Technology, Ryukoku University,\\
\qquad\quad
1-5 Yokotani, Seta Oe-cho, Otsu, Shiga 520-2194, Japan.\\
\noindent 
E-mail: {\tt kawakami@math.ryukoku.ac.jp}\\

\newpage
\section{Introduction}
This paper is concerned with the large time behavior of a solution 
to the Cauchy problem for an inhomogeneous fractional diffusion equation
\begin{equation}
\label{eq:1.1}
	\partial_t u+(-\Delta)^{\frac{\theta}{2}}u=f(x,t)\quad\mbox{in}\quad{\mathbb R}^N\times(0,\infty),
	\quad
	u(x,0)=\varphi(x)\quad\mbox{in}\quad{\mathbb R}^N,
\end{equation}
where $N\ge 1$, $\partial_t:=\partial/\partial t$, $0<\theta<2$ and 
$$
\varphi\in L_K^1:=L^1({\mathbb R}^N,(1+|x|)^K\,dx)\quad\mbox{with}\quad K\ge 0. 
$$
Here $(-\Delta)^{\theta/2}$ is the fractional power of the Laplace operator.  
Inhomogeneous fractional diffusion equation~\eqref{eq:1.1} 
appears in the study of various nonlinear problems with anomalous diffusion, 
the Laplace equation with a dynamical boundary condition, and so on. 
Under suitable integrability conditions on the inhomogeneous term $f$, 
the solution~$u$ to problem~\eqref{eq:1.1} behaves like a suitable multiple of the fundamental solution~$G_\theta$
to the linear fractional diffusion equation
$$
	\partial_t v+(-\Delta)^{\frac{\theta}{2}}v=0\quad\mbox{in}\quad{\mathbb R}^N\times(0,\infty)
$$
as $t\to\infty$. 
In this paper we obtain the higher order asymptotic expansions (HOAE) of the large time behavior of 
the solution~$u$. Furthermore, we study the precise description of the large time behavior of solutions to the Cauchy problem 
for nonlinear fractional diffusion equations such as 
\begin{equation}
\label{eq:1.2}
	\partial_t u+(-\Delta)^{\frac{\theta}{2}}u=\lambda|u|^{p-1}u
	\quad\mbox{in}\quad{\mathbb R}^N\times(0,\infty),
	\quad
	u(x,0)=\varphi(x)\quad\mbox{in}\quad{\mathbb R}^N,
\end{equation}
where $\lambda\in{\mathbb R}$, $p>1$ and $\varphi\in L_K^1$ with $K\ge 0$. 
This paper is an improvement of \cite{IKM} and 
it corresponds a fractional version of the papers~\cites{IIK, IK, IKK02}.
\vspace{5pt}

The large time behavior of solutions to nonlinear parabolic equations has been studied extensively 
in many papers by various methods. 
Here we just refer to the papers 
\cites{AJY, BK, FK, FIK, FM, HKN, IIK, IK, IKK01, IKK02, IKM, Iwa, YM, YS01, YS02, YS03, YT1, YT2}, 
which are closely related to this paper. 
Among others, 
in \cites{IIK, IK, IKK02}, 
HOAE of solutions behaving like suitable multiples of the Gauss kernel 
have already been well established.  The property that 
$$
	\bigcup_{t>0}e^{t\Delta}L_K^1\subset L_K^1
	\quad\mbox{for}\quad K\ge 0
$$
plays an important role in \cites{IIK, IK, IKK02} and it follows 
from the exponential decay of the Gauss kernel at the space infinity. 
For fractional diffusion equations, 
if $0\le K<\theta$, then
\begin{equation}
\label{eq:1.3}
	\bigcup_{t>0}e^{-t(-\Delta)^{\theta/2}}L_K^1\subset L_K^1
\end{equation}
holds and the arguments in \cites{IIK, IK, IKK02} are also applicable to fractional diffusion equations. 
However, if $K\ge\theta$, then property~\eqref{eq:1.3} fails. 
This fact prevents to establish analogous asymptotic expansions of solutions to the case of $\theta=2$. 
In \cite{IKM} the authors of this paper and Michihisa studied a mechanism for property~\eqref{eq:1.3} to fail in the case of $K\ge\theta$, 
and obtained HOAE of $e^{-(-\Delta)^{\theta/2}}\varphi$. 
This argument is applicable to 
the study of HOAE of solutions 
to inhomogeneous fractional diffusion equations and nonlinear fractional diffusion equations,  
however HOAE of \cite{IKM} to problem~\eqref{eq:1.1} do not have refined forms.

In this paper we improve and refine arguments in \cite{IKM} by 
taking into an account of the Taylor expansion of the kernel $G_\theta$ with respect to 
both of the space and the time variables, and obtain 
HOAE of solutions to inhomogeneous fractional diffusion equations and nonlinear fractional  parabolic equations. 
Our arguments also reveal a mechanism 
for the solution $u$ to problem~\eqref{eq:1.1} to break the property that
$u(t)\in L^1_K$ for $t>0$. 

We introduce some notations. 
Set $\mathbb N_0:={\mathbb N}\cup\{0\}$.
For any $k\ge 0$, let $[k]\in\mathbb N_0={\mathbb N}\cup\{0\}$
be such that $k-1<[k]\le k$. 
Let $\nabla:=(\partial/\partial x_1,\dots,\partial/\partial x_N)$. 
For any multi-index $\alpha\in{\mathbb M}:=\mathbb N_0^N$, 
set
$$
	|\alpha|:=\displaystyle{\sum_{i=1}^N}\alpha_i,\quad
	\alpha!:=\prod_{i=1}^N\alpha_i!,\quad
	x^\alpha:=\prod_{i=1}^N x_i^{\alpha_i},\quad
	\partial_x^\alpha:=
	\frac{\partial^{|\alpha|}}{\partial x_1^{\alpha_1}\cdots\partial x_N^{\alpha_N}}.
$$
For any $\alpha=(\alpha_1,\dots,\alpha_N)$, $\beta=(\beta_1,\dots,\beta_N)\in{\mathbb M}$, 
we say 
$\alpha\le\beta$
if $\alpha_i\le\beta_i$ for all $i\in\{1,\dots,N\}$. 
Let $1\le q\le \infty$ and $K\ge 0$. 
Let $\|\cdot\|_q$ be the usual norm of $L^q:=L^q({\mathbb R}^N)$. 
Set 
$$
	|||f|||_{q,K}:=\|f_K\|_q\quad\mbox{with}\quad f_K(x):=|x|^K f(x).
$$
Let 
$$
	f\in L_K^q:=\left\{f\in L^q\,:\,\|f\|_{L^q_K}<\infty\right\},
	\quad
	\mbox{where}\quad
	\|f\|_{L^q_K}:=\|f\|_q+|||f|||_{q,K}.
$$
For any $f\in L^1_K$ and $\alpha\in\mathbb M$ with $|\alpha|\le K$, 
set 
$$
	M_\alpha(f):=\int_{\mathbb R^N}x^\alpha f(x)\,dx.
$$

We are ready to state our main results on the asymptotic expansions of solutions to inhomogeneous fractional diffusion equations. 
In what follows, set $K_\theta:=[K/\theta]$. Furthermore, set 
$$
g_{\alpha,m}(x,t):=\frac{(-1)^{|\alpha|+m}}{\alpha!m!}(\partial_t^m\partial_x^\alpha G_\theta)(x,t+1)
$$
for $(x,t)\in{\mathbb R}^N\times(0,\infty)$, 
where $\alpha\in{\mathbb M}$ and $m\in{\mathbb N}_0$. 
\begin{theorem}
\label{Theorem:1.1}
	Let $N\ge1$, $0<\theta<2$, $0\le\ell\le K$, and $1\le q\le\infty$.
	Let $\varphi\in L^1_K$ and $f$ be a measurable function in $\mathbb R^N\times(0,\infty)$ such that 
	\begin{equation}
	\label{eq:1.4}
	E_{K,q}[f]\in L^1_{{\rm loc}}(0,\infty),
	\end{equation}
	where
	\begin{equation}
	\label{eq:1.5}
		\begin{aligned}
		E_{K,q}[f](t)
	 	& := (t+1)^{\frac{K}{\theta}}
		\left[t^{\frac{N}{\theta}\left(1-\frac{1}{q}\right)}\|f(t)\|_q+\|f(t)\|_1\right]\\
	 	& \qquad
		+t^{\frac{N}{\theta}\left(1-\frac{1}{q}\right)}|||f(t)|||_{q,K}+|||f(t)|||_{1,K}
		\quad\mbox{for}\quad t>0. 
		\end{aligned}
	\end{equation}
	Let $u\in C({\mathbb R}^N\times(0,\infty))$ be a solution to problem~\eqref{eq:1.1}, that is, 
	$u$ satisfies
	$$
		u(x,t)
		=\int_{{\mathbb R}^N}G_\theta(x-y,t)\varphi(y)\,dy
		+\int_0^t\int_{{\mathbb R}^N}G_\theta(x-y,t-s)f(y,s)\,dy\,ds
	$$
	for $(x,t)\in{\mathbb R}^N\times(0,\infty)$. 
	Then
	\begin{equation}
	\label{eq:1.6}
 		\sup_{0<t<\tau}\,t^{\frac{N}{\theta}\left(1-\frac{1}{q}\right)}|||u(t)-w(t)|||_{q,\ell}<\infty
		\quad\mbox{for}\quad\tau>0,
	\end{equation}
	where 
	\begin{equation}
	\label{eq:1.7}
		w(x,t):=
		\sum_{m=0}^{K_\theta}\sum_{|\alpha|\le K}
		\bigg\{M_\alpha(\varphi)
		+\int_0^t(s+1)^m M_\alpha(f(s))\,ds\bigg\}\,g_{\alpha,m}(x,t).
	\end{equation}
	Furthermore, there exists $C>0$ such that, 
	for any $\varepsilon>0$ and $T>0$, 
	\begin{equation}
	\label{eq:1.8}
 	 	t^{\frac{N}{\theta}\left(1-\frac{1}{q}\right)-\frac{\ell}{\theta}}|||u(t)-w(t)|||_{q,\ell}\\
  	 	\le\varepsilon t^{-\frac{K}{\theta}}
		+Ct^{-\frac{K}{\theta}}\int_T^tE_{K,q}[f](s)\,ds
	\end{equation}
	holds for large enough $t>0$.
	In particular, if 
	$$
		\int_0^\infty E_{K,q}[f](s)\,ds<\infty,
	$$
	then 
	\begin{equation}
	\label{eq:1.9}
		\lim_{t\to\infty}t^{\frac{N}{\theta}\left(1-\frac{1}{q}\right)+\frac{K-\ell}{\theta}}|||u(t)-w(t)|||_{q,\ell}=0.
	\end{equation}
\end{theorem}
Theorem~\ref{Theorem:1.1} corresponds to \cite{IKK02}*{Theorems~1.1, 1.2} for $\theta=2$ 
and it is an improvement of \cite{IKM}*{Theorem~3.1~(ii)}. 
Our asymptotic profile $w$ has a pretty simpler form than that of \cite{IKM}. 
(See Remarks~\ref{Remark:3.1} and \ref{Remark:5.2}.) 
We also remark that, under condition~\eqref{eq:1.4}, 
both of $u(\cdot,t)$ and $w(\cdot,t)$ do not necessarily belong to $L^q_\ell$, 
while $u(t)-w(t)\in L^q_\ell$. 
In other words, the function~$w$ may break the property that $u(t)\in L^q_\ell$ for $t>0$. 
Furthermore, we have: 
\begin{corollary}
\label{Corollary:1.1}
Assume the same conditions as in Theorem~{\rm\ref{Theorem:1.1}}. 
Let $u\in C({\mathbb R}^N\times(0,\infty))$ be a solution to problem~\eqref{eq:1.1}. 
Then there exists $R>0$ such that 
$$
	u(t)\in \left\{h+\sum_{(\alpha,m)\in\Lambda_K^q}a_{\alpha,m}g_{\alpha,m}(x,t)\,:\,
	\mbox{$h\in L^q_K$ with $\|h\|_{L^q_K}\le R$},\,\{a_{\alpha,m}\}\subset[-R,R]\right\}
$$
for $t>0$, where 
$\Lambda_K^q:=\{(\alpha,m)\in{\mathbb M}\times{\mathbb N}_0\,:\,g_{\alpha,m}(\cdot,0)\not\in L^q_K\}$. 
\end{corollary}

We explain the idea of the proof of Theorem~\ref{Theorem:1.1}. 
We improve and refine arguments in the previous papers~\cites{IIK, IK, IKK02, IKM} 
to obtain HOEA of the solution~$u$ to problem~\eqref{eq:1.1}, in  particular, the integral term
$$
\int_0^t\int_{{\mathbb R}^N}G_\theta(x-y,t-s)f(y,s)\,dy\,ds.
$$
In \cite{IKM}, 
following the arguments in \cites{IIK, IK, IKK02}, 
the authors of this paper and Michihisa expanded 
the integral kernel $G_\theta(x-y,t-s)$ by the Taylor expansions with respect to the space derivatives of $G_\theta(x,t-s)$. 
Then the slow decay of $G_\theta(x,t)$ 
makes difficult to obtain refined HOAE of the solution~$u$ to problem~\eqref{eq:1.1}. 
In this paper we expand the integral kernel $G_\theta(x-y,t-s)$ 
by the Taylor expansions with respect to both of space and time variable derivatives of $G_\theta(x,t)$. 
(This is the same sprit as in \cite{FM}.) 
Indeed, we introduce the following integral kernels by the use of the Taylor expansions of $G_\theta$:
\begin{equation}
\label{eq:1.10}
	\begin{aligned}
 	{\mathcal S}^m_\ell(x,y,t) & :=(\partial_t^m G_\theta)(x-y,t) 
  	-\sum_{|\alpha| \le \ell}\frac{(-1)^{|\alpha|}}{\alpha !} 
  	(\partial_t^m\partial^\alpha_x G_\theta)(x,t)y^\alpha\\
 	& \,\, =\frac{1}{[\ell]!}\int_0^1(1-\tau)^{[\ell]}\frac{\partial^{[\ell]+1}}{\partial\tau^{[\ell]+1}}
	(\partial_t^m G_\theta)(x-\tau y,t)\,d\tau,\\
 	{\mathcal T}(x,y,t,s) & :=G_\theta(x-y,t-s)-\sum_{m=0}^{K_\theta} \frac{(-1)^{m}}{m!} 
  	(\partial_t^m G_\theta)(x-y,t)s^m\\
 	& \,\,=\frac{1}{K_\theta!}
	\int_0^1(1-\tau)^{K_\theta}\frac{\partial^{K_\theta+1}}{\partial\tau^{K_\theta+1}}G_\theta(x-y,t-\tau s)\,d\tau,
 	\end{aligned}
\end{equation}
for $x$, $y\in{\mathbb R}^N$ and $0\le s<t$, where $0\le\ell\le K$ and $m\in{\mathbb N}_0$. 
Then
$$
	\begin{aligned}
 	 & G_\theta(x-y,t-s)\\
	 & =\sum_{m=0}^{K_\theta} \frac{(-1)^{m}}{m!}(\partial_t^m G_\theta)(x-y,t)s^m+{\mathcal T}(x,y,t,s)\\
 	& =\sum_{m=0}^{K_\theta} \sum_{|\alpha| \le \ell}\frac{(-1)^{|\alpha|+m}}{\alpha ! m!} 
  	(\partial_t^m\partial^\alpha_x G_\theta)(x,t)y^\alpha s^m
	+\sum_{m=0}^{K_\theta} \frac{(-1)^m}{m!}{\mathcal S}_\ell^m(x,y,t) s^m+{\mathcal T}(x,y,t,s).
	\end{aligned}
$$
Furthermore, 
\begin{equation}
\label{eq:1.11}
\begin{split}
	{\mathcal R}(x,y,t,s) & :={\mathcal T}(x,y,t,s)+\sum_{m=0}^{K_\theta} \frac{(-1)^m}{m!}{\mathcal S}_K^m(x,y,t) s^m\\
	 & \,=G_\theta(x-y,t-s) 
	-\sum_{m=0}^{K_\theta}\sum_{|\alpha| \le K}\frac{(-1)^{|\alpha|+m}}{\alpha ! m!} 
	(\partial_t^m\partial^\alpha_x G_\theta)(x,t)y^\alpha s^m.
\end{split}
\end{equation}
Then it follows from \eqref{eq:1.7} that
\begin{equation}
\label{eq:1.12}
	\begin{aligned}
 	u(x,t)-w(x,t)
	 & =\int_{{\mathbb R}^N}{\mathcal R}(x,y,t+1,1)\varphi(y)\,dy\\
	& +\int_0^t\int_{{\mathbb R}^N} {\mathcal R}(x,y,t+1,s+1)f(y,s)\,dy\,ds.
	\end{aligned}
\end{equation}
Thanks to the decay of the derivatives of $G_\theta$ and \eqref{eq:1.10}, 
we see that 
\begin{equation*}
	\begin{aligned}
 	& {\mathcal T}(x,y,t,s)=O(|x-y|^{-N-K-\varepsilon})\quad\mbox{as}\quad |x-y|\to\infty,\\
 	& {\mathcal S}^m_K(x,y,t)=O(|x|^{-N-K-\varepsilon})\quad\mbox{as}\quad |x|\to\infty,
	\end{aligned}
\end{equation*}
for some $\varepsilon>0$. 
These decay of the integral kernels at the space infinity enables us 
to establish HOAE of solutions to problem~\eqref{eq:1.1} 
and to obtain Theorem~\ref{Theorem:1.1}. 
These arguments require delicate integral estimates on the integral kernels 
${\mathcal S}_\ell^m$ and ${\mathcal T}$. 
\vspace{5pt}

Theorem~\ref{Theorem:1.1} is applicable to problem~\eqref{eq:1.2} 
and it gives asymptotic profiles of solutions to problem~\eqref{eq:1.2} 
as a linear combination of the derivatives of $G_\theta$ (see Theorem~\ref{Theorem:5.1}).
Furthermore, 
taking a suitable approximation of the nonlinear term in problem~\eqref{eq:1.2}, 
we obtain refined asymptotic expansions of the solution to problem~\eqref{eq:1.1}
(see Theorem~\ref{Theorem:5.2}).
Here we state the following result, which is a variation of Theorem~\ref{Theorem:5.2}.
\begin{theorem}
\label{Theorem:1.2}
	Let $N\ge 1$, $0<\theta<2$, $\lambda\in{\mathbb R}$, and $\varphi\in L^1_K\cap L^\infty$ with $K\ge0$.
	Let $u\in C({\mathbb R}^N\times(0,\infty))$ be a solution to problem~\eqref{eq:1.2} with $p>1+\theta/N$
	and satisfy
	\begin{equation}
	\label{eq:1.13}
		\sup_{t>0}\,(t+1)^{\frac{N}{\theta}}\|u(t)\|_\infty<\infty.
	\end{equation}
	Then there exists $M_*\in{\mathbb R}$ such that 
  		$$
  			M_*:=\lim_{t\to\infty}\int_{{\mathbb R}^N}u(x,t)\,dx
  			=\int_{{\mathbb R}^N}\varphi(x)\,dx+\int_0^\infty\int_{{\mathbb R}^N}F(u(x,t))\,dx\,dt,
  		$$ 
		where $F(u(x,t)):=\lambda|u(x,t)|^{p-1}u(x,t)$.
		
	Assume $N(p+\theta)>N+K$ and $\varphi\in L^\infty_k$ with $k=\min\{N+\theta,K\}$. 
	Let $1\le q\le\infty$. 
	Then 
	$$
		\sup_{t>0}\,(t+1)^{\frac{N}{\theta}\left(1-\frac{1}{q}\right)-\frac{\ell}{\theta}}|||u(t)|||_{q,\ell}<\infty,
	$$
	where $0\le \ell \le K$ with $0\le \ell< \theta+N(1-1/q)$. 
	Furthermore, for any $\sigma>0$ 
	$$
		\begin{aligned}
  		 &\sup_{t>0}\,t^{\frac{N}{\theta}\left(1-\frac{1}{q}\right)-\frac{\ell}{\theta}}
		|||u(t)-v(t)|||_{q,\ell}<\infty,
 		\\
		&
  		t^{\frac{N}{\theta}\left(1-\frac{1}{q}\right)-\frac{\ell}{\theta}}|||u(t)-v(t)|||_{q,\ell}
		=o\left(t^{-\frac{K}{\theta}}\right)
  		+O\left(t^{-\frac{K}{\theta}}\int_1^t s^{\frac{K}{\theta}-A_p}h_\sigma(s)\,ds\right)
		\quad\mbox{as}\quad t\to\infty,
		\end{aligned}
 	$$
  	where $0\le\ell\le K$. 
	Here
  	$$
  		\begin{aligned}
		 & v(x,t):=
		\sum_{m=0}^{K_\theta}\sum_{|\alpha|\le K}c_{\alpha,m}(t)g_{\alpha,m}(x,t)
   		+\int_0^t e^{-(t-s)(-\Delta)^{\frac{\theta}{2}}}F_\infty(s)\,ds,\\
		 & c_{\alpha,m}(t):=M_\alpha(\varphi)
		+\int_0^t (s+1)^m M_\alpha(F(u(s))-F_\infty(s))\,ds,\\
		 & F_\infty(x,t):=F\left(M_*G_\theta(x,t+1)\right),\quad
  		h_\sigma(t):=t^{-(A_p-1)+\sigma}+t^{-1}+t^{-\frac{1}{\theta}}.
  		\end{aligned}
  	$$	
\end{theorem}
Theorem~\ref{Theorem:1.2} corresponds to \cite{IK}*{Corollary~1.1} for $\theta=2$. 
See Remark~\ref{Remark:5.1} for condition~\eqref{eq:1.13}. 
\vspace{3pt}

The rest of this paper is organized as follows. 
In Section~2 we collect some properties of the fundamental solution $G_\theta$. 
In Section~3 we obtain some estimates on the integral kernels ${\mathcal S}_\ell^m(x,y,t)$ 
and ${\mathcal T}(x,y,t,s)$, 
and prove Theorem~\ref{Theorem:1.1} and Corollary~\ref{Corollary:1.1}.
In Section~4 we apply Theorem~\ref{Theorem:1.1} 
to obtain HOAE of solutions to 
the Cauchy problem for a convection type inhomogeneous fractional diffusion equation.
In Section~5 we apply Theorem~\ref{Theorem:1.1} to 
study HOAE of solutions to the Cauchy problem 
for nonlinear fractional diffusion equations. 
Furthermore, we prove Theorem~\ref{Theorem:1.2}.
\section{Preliminaries}
We recall some properties of the fundamental solution $G_\theta=G_\theta(x,t)$,
In what follows,
by the letter $C$
we denote generic positive constants (independent of $x$ and $t$)
and they may have different values also within the same line.

Let $0<\theta<2$.
The fundamental solution $G_\theta=G_\theta(x,t)$ is represented by 
$$
	G_\theta(x,t)=(2\pi)^{-\frac{N}{2}}\int_{{\mathbb R}^N}e^{ix\cdot \xi}e^{-t|\xi|^\theta}\,d\xi,
	\quad
	(x,t)\in{\mathbb R}^N\times(0,\infty).
$$
Then we have: 
\begin{itemize}
\item[({\bf G})]
  	$G_\theta=G_\theta(x,t)$ is a positive smooth function in ${\mathbb R}^N\times(0,\infty)$ with the following properties:
  	\begin{itemize}
  	\item[(i)] 
		$\displaystyle{G_\theta(x,t)=t^{-\frac{N}{\theta}}G_\theta(t^{-\frac{1}{\theta}}x,1)}$ 
		for $x\in{\mathbb R}^N$ and $t>0$;
	\item[(ii)]
		$\displaystyle{\sup_{x\in{\mathbb R}^N}\,(1+|x|)^{N+\theta+|\alpha|}
		|(\partial_x^\alpha G_\theta)(x,1)|<\infty}$ 
		for $\alpha\in{\mathbb M}$;
	\item[(iii)]
		$G_\theta(\cdot,1)$ is radially symmetric and decreasing with respect to $r:=|x|$. Furthermore, 
		$$
		\liminf_{|x|\to+\infty}\,(1+|x|)^{N+\theta+j}(\partial_r^j G_\theta)(x,1)>0,\quad j\in\mathbb N_0;
		$$
	\item[(iv)] 
		$\displaystyle{G_\theta(x,t)=\int_{{\mathbb R}^N}G_\theta(x-y,t-s)\,G_\theta(y,s)\,dy}$ 
		for $x\in\mathbb R^N$ and $t>s>0$;
	\item[(v)]
		$\displaystyle{\int_{\mathbb R^N}G_\theta(x,t)\,dx=1}$ for $t>0$.
  	\end{itemize}
\end{itemize}
See \cites{BT, BK}. (See also \cites{IKK01, IKK02, S}.) 

Let $\alpha\in{\mathbb M}$ and $m\in{\mathbb N}_0$. 
Let 
$$
	1\le q\le\infty,\quad
	0\le \ell<\theta m'+|\alpha|+N\left(1-\frac{1}{q}\right)
	\quad\mbox{with}\quad m':=\max\{m,1\}.
$$
It follows from ({\bf G})-(i), (ii) and \cite{IKM}*{Lemma~2.1} that
\begin{equation}
\label{eq:2.1}
	|(\partial_t^m\partial_x^\alpha G_\theta)(x,t)|
	\le Ct^{-\frac{N+|\alpha|}{\theta}-m}
	\big(1+t^{-\frac{1}{\theta}}|x|\big)^{-(N+\theta m'+|\alpha|)}
\end{equation}
for $x\in{\mathbb R}^N$ and $t>0$.
This implies that
\begin{equation}
\label{eq:2.2}	
	\sup_{t>0}\,
	t^{\frac{N}{\theta}\left(1-\frac{1}{q}\right)+\frac{|\alpha|-\ell}{\theta}+m}
	|||(\partial_t^m\partial_x^\alpha G_\theta)(t)|||_{q,\ell}<\infty.
\end{equation}
\begin{lemma}
\label{Lemma:2.1}
Let $1\le q\le r\le\infty$, $\alpha\in{\mathbb M}$, and $m\in{\mathbb N}_0$. 
Let 
\begin{equation}
\label{eq:2.3}
0\le \ell<\theta m'+|\alpha|+N\left(\frac{1}{q}-\frac{1}{r}\right).
\end{equation}
Then there exists $C>0$ such that 
$$
	t^{\frac{N}{\theta}\left(\frac{1}{q}-\frac{1}{r}\right)+\frac{|\alpha|}{\theta}+m}
	\biggr|\biggr|\biggr|\,\partial_t^m\partial_x^\alpha e^{-t(-\Delta)^{\theta/2}}\varphi\,
	\biggr|\biggr|\biggr|_{r,\ell}
	\le Ct^{\frac{\ell}{\theta}}\|\varphi\|_q+C|||\varphi|||_{q,\ell}
$$
for $\varphi\in L^q_\ell$ and $t>0$. Here 
$$
	\left[e^{-t(-\Delta)^{\theta/2}}\varphi\right](x):=\int_{{\mathbb R}^N}G_\theta(x-y,t)\varphi(y)\,dy,
	\quad (x,t)\in{\mathbb R}^N\times(0,\infty).
$$
\end{lemma}
{\bf Proof.}
Assume \eqref{eq:2.3}. It follows that 
$$
	|x|^\ell\left[\partial_t^m\partial_x^\alpha e^{-t(-\Delta)^{\theta/2}}\varphi\right](x)
	\le C\int_{{\mathbb R}^N}\left[|x-y|^\ell+|y|^\ell\right]
	|(\partial_t^m\partial_x^\alpha G_\theta)(x-y,t)||\varphi(y)|\,dy
$$
for $(x,t)\in{\mathbb R}^N\times(0,\infty)$. 
The Young inequality together with \eqref{eq:2.2} implies that
$$
	\begin{aligned}
 	& \biggr|\biggr|\biggr|\,\partial_t^m\partial_x^\alpha e^{-t(-\Delta)^{\theta/2}}\varphi\,	
	\biggr|\biggr|\biggr|_{r,\ell}\\
 	& \le C|||\partial_t^m\partial_x^\alpha G_\theta(t)|||_{p,\ell}\|\varphi\|_q
	+C\|	\partial_t^m\partial_x^\alpha G_\theta(t)\|_p|||\varphi|||_{q,\ell}\\
	 & \le Ct^{-\frac{N}{\theta}\left(\frac{1}{q}-\frac{1}{r}\right)-\frac{|\alpha|}{\theta}-m+\frac{\ell}{\theta}}
	 \|\varphi\|_q
 	+Ct^{-\frac{N}{\theta}\left(\frac{1}{q}-\frac{1}{r}\right)-\frac{|\alpha|}{\theta}-m}|||\varphi|||_{q,\ell}
	\end{aligned}
$$
for $t>0$, where $p\in[1,\infty]$ with $1/r=1/p+1/q-1$.
Then we obtain the desired inequality, and the proof is complete. 
$\Box$
%
\section{Proof of Theorem~\ref{Theorem:1.1}}
In this section we prove Theorem~\ref{Theorem:1.1}. 
We first prepare the following lemma.
\begin{lemma}
\label{Lemma:3.1}
	Assume the same conditions as in Theorem~{\rm\ref{Theorem:1.1}}. 
	Then 
	$$
		t^{\frac{N}{\theta}\left(1-\frac{1}{r}\right)}(t+1)^{\frac{K-\ell}{\theta}}|||f(t)|||_{r,\ell}\le E_{K,q}[f](t),
		\quad t>0,
	$$
	where $0\le\ell\le K$ and $1\le r\le q$.
\end{lemma}
{\bf Proof.}
Let $0\le\ell\le K$ and $1\le q\le\infty$. 
It follows that 
$$
	(t+1)^{-\frac{\ell}{\theta}}|x|^\ell\le C+C(t+1)^{-\frac{K}{\theta}}|x|^K,
	\quad (x,t)\in{\mathbb R}^N\times(0,\infty).
$$
This together with \eqref{eq:1.5} implies that
$$
	\begin{aligned}
 	(t+1)^{-\frac{\ell}{\theta}}|||f(t)|||_{r,\ell}
 	& \le C\|f(t)\|_r+C(t+1)^{-\frac{K}{\theta}}|||f(t)|||_{r,K}\\
 	& \le C\|f(t)\|_1^\lambda\|f(t)\|_q^{1-\lambda}
	+C(t+1)^{-\frac{K}{\theta}}|||f(t)|||_{1,K}^\lambda|||f(t)|||_{q,K}^{1-\lambda}\\
 	& \le Ct^{-\frac{N}{\theta}\left(1-\frac{1}{r}\right)}(t+1)^{-\frac{K}{\theta}}E_{K,q}[f](t),
	\quad t>0,
 	\end{aligned}
$$
where $1/r=\lambda+(1-\lambda)/q$. 
Thus Lemma~\ref{Lemma:3.1} follows. 
$\Box$
\vspace{5pt}

Next we prove a lemma on the integral kernel ${\mathcal S}^m_\ell(x,y,t)$.
\begin{lemma}
\label{Lemma:3.2}
	Let $m\in{\mathbb N_0}$, $0\le\ell \le K$, $1\le q\le\infty$, and $j=0,1$. 
	\begin{itemize}
	\item[{\rm (a)}]
		There exists $C_1>0$ such that
		$$
			|||\nabla^j{\mathcal S}^m_\ell(\cdot,y,t)|||_{q,\ell}
			\le C_1t^{-\frac{N}{\theta}\left(1-\frac{1}{q}\right)-m-\frac{j}{\theta}}|y|^\ell,
			\quad (y,t)\in{\mathbb R}^N\times(0,\infty).
		$$
	\item[{\rm (b)}] 
		There exists $C_2>0$ such that 
		$$
			\left|\left|\left|\,
			\int_{\mathbb R^N}\nabla^j{\mathcal S}^m_K(\cdot,y,t)\varphi(y)\,dy\,\right|\right|\right|_{q,\ell}
			\le C_2t^{-\frac{N}{\theta}\left(1-\frac{1}{q}\right)-m-\frac{K+j-\ell}{\theta}}
			|||\varphi|||_{1,K},\quad t>0,
		$$
		for $\varphi\in L^1_K$. 
	\item[{\rm (c)}] 
		Let $\varphi\in L^1_K$. Then 
		$$
			\lim_{t\to\infty}t^{\frac{N}{\theta}\left(1-\frac{1}{q}\right)+m+\frac{K+j-\ell}{\theta}}
			\left|\left|\left|\,
			\int_{\mathbb R^N}\nabla^j{\mathcal S}^m_K(\cdot,y,t)\varphi(y)\,dy\,
			\right|\right|\right|_{q,\ell}=0.
		$$
	\end{itemize}
\end{lemma}
{\bf Proof.} 
Let $0\le\ell\le K$, $1\le q\le\infty$, and $j=0,1$. 
We prove assertion~(a). 
Let $x$, $y\in{\mathbb R}^N$ and $t>0$. 
It follows that 
$$
	|x-\tau y| \ge |x|-|y| \ge |x|/2
	\quad\mbox{if}\quad
	|x| \ge 2|y|\quad\mbox{and}\quad
	0 \le \tau \le 1.
$$
Then, by \eqref{eq:1.10} and \eqref{eq:2.1} we have
\begin{equation*}
	\begin{aligned}
	&
	|x|^\ell |\nabla^j{\mathcal S}^m_\ell(x,y,t)| \\
 	& \le C\int_0^1|x|^\ell 
	\left|(\partial_t^m\nabla^{[\ell]+j+1} G_\theta)(x-\tau y,t) \right| |y|^{[\ell]+1}\,d\tau \\
 	& \le C|y|^\ell \int_0^1|x|^{[\ell]+1} \ t^{-\frac{N}{\theta}-\frac{[\ell]+j+1}{\theta}-m} 
	\left(1+t^{-\frac{1}{\theta}}|x-\tau y| \right)^{-(N+\theta m'+[\ell]+j+1)}\,d\tau\\
 	& \le C|y|^\ell (t^{-\frac{1}{\theta}}|x|)^{[\ell]+1} \ t^{-\frac{N}{\theta}-m-\frac{j}{\theta}} 
	\left(1+t^{-\frac{1}{\theta}}\frac{|x|}{2} \right)^{-(N+\theta m'+[\ell]+j+1)}\\
 	& \le C|y|^\ell t^{-\frac{N}{\theta}-m-\frac{j}{\theta}} 
	\left(1+t^{-\frac{1}{\theta}}\frac{|x|}{2} \right)^{-(N+\theta m'+j)}
	\end{aligned}
\end{equation*}
if $|x|\ge 2|y|$. 
Similarly, by \eqref{eq:1.10} we have
\begin{equation*}
	\begin{aligned}
	|x|^\ell |\nabla^j{\mathcal S}^m_\ell(x,y,t)|
 	& \le  |x|^\ell |(\partial_t^m \nabla^jG_\theta)(x-y,t)|
	+C\sum_{|\alpha|\le\ell}|x|^\ell |(\partial_t^m\partial_x^\alpha \nabla^jG_\theta)(x,t)||y|^{|\alpha|}\\
 	& \le (2|y|)^\ell |(\partial_t^m \nabla^jG_\theta)(x-y,t)|
	+C\sum_{|\alpha|\le\ell}  |y|^\ell |x|^{|\alpha|} 
 	|(\partial_t^m\partial_x^\alpha \nabla^jG_\theta)(x,t)|
	\end{aligned}
\end{equation*}
if $|x|<2|y|$. 
These together with \eqref{eq:2.2} imply that 
$$
	|||\nabla^j{\mathcal S}^m_\ell(\cdot,y,t)|||_{q,\ell}
	\le Ct^{-\frac{N}{\theta}\left(1-\frac{1}{q}\right)-m-\frac{j}{\theta}}|y|^\ell.
$$
Thus assertion~(a) follows. 

We prove assertions~(b) and (c). 
Let $\varphi\in L^1_K$, $0\le\ell\le K$, and $R>0$. 
It follows from \eqref{eq:1.10} that
\begin{equation*}
	\begin{aligned}
 	& \biggr|\biggr|\biggr|\,\int_{\{|y|\ge R^{\frac{1}{\theta}}\}}
 	|\nabla^j{\mathcal S}^m_K(\cdot,y,t)||\varphi(y)|\,dy\,\biggr|\biggr|\biggr|_{q,\ell}\\
	& \le\,\biggr|\biggr|\biggr|\,\int_{\{|y|\ge R^{\frac{1}{\theta}}\}}
 	|\nabla^j{\mathcal S}^m_\ell(\cdot,y,t)||\varphi(y)|\,dy\,\biggr|\biggr|\biggr|_{q,\ell}\\
	 & \qquad
	 +C\sum_{\ell<|\alpha|\le K}\biggr|\biggr|\biggr|\,\int_{\{|y|\ge R^{\frac{1}{\theta}}\}}
 	|(\partial_t^m\partial_x^\alpha \nabla^jG_\theta)(\cdot,t)||y|^{|\alpha|}|\varphi(y)|\,dy\,\biggr|\biggr|\biggr|_{q,\ell}\\
	 & \le \int_{\{|y|\ge R^{\frac{1}{\theta}}\}}
 	|||\nabla^j{\mathcal S}^m_\ell(\cdot,y,t)|||_{q,\ell}|\varphi(y)|\,dy\\
	 & \qquad
	 +C\sum_{\ell<|\alpha|\le K}\int_{\{|y|\ge R^{\frac{1}{\theta}}\}}
 	|||(\partial_t^m\partial_x^\alpha \nabla^jG_\theta)(\cdot,t)|||_{q,\ell}|y|^{|\alpha|}|\varphi(y)|\,dy.
	\end{aligned}
\end{equation*}
This together with \eqref{eq:2.2} and assertion~(a) implies that
\begin{equation}
\label{eq:3.1}
	\begin{aligned}
 	& t^{\frac{N}{\theta}\left(1-\frac{1}{q}\right)+m+\frac{j}{\theta}}\biggr|\biggr|\biggr|\,\int_{\{|y|\ge R^{\frac{1}{\theta}}\}}
 	|\nabla^j{\mathcal S}^m_K(\cdot,y,t)||\varphi(y)|\,dy\,\biggr|\biggr|\biggr|_{q,\ell}\\
 	& \le C\bigg(\int_{\{|y|\ge R^{\frac{1}{\theta}}\}}|y|^\ell|\varphi(y)|\,dy
 	+\sum_{\ell<|\alpha|\le K}t^{-\frac{|\alpha|-\ell}{\theta}}
	\int_{\{|y|\ge R^{\frac{1}{\theta}}\}}|y|^{|\alpha|}|\varphi(y)|\,dy\bigg)\\
 	& \le C\bigg(\int_{\{|y|\ge R^{\frac{1}{\theta}}\}}|y|^\ell
	\left(\frac{|y|}{R^{\frac{1}{\theta}}}\right)^{K-\ell}|\varphi(y)|\,dy\\
 	& \qquad\qquad\qquad\qquad\qquad
 	+\sum_{\ell<|\alpha|\le K}t^{-\frac{|\alpha|-\ell}{\theta}}
 	\int_{\{|y|\ge R^{\frac{1}{\theta}}\}}
	\left(\frac{|y|}{R^{\frac{1}{\theta}}}\right)^{K-|\alpha|}|y|^{|\alpha|}|\varphi(y)|\,dy\bigg)\\
 	& =Ct^{-\frac{K-\ell}{\theta}}
 	\biggr((R^{-1}t)^{\frac{K-\ell}{\theta}}
	+\sum_{\ell<|\alpha|\le K}(R^{-1}t)^{\frac{K-|\alpha|}{\theta}}\biggr)
 	\int_{\{|y|\ge R^{\frac{1}{\theta}}\}}|y|^K|\varphi(y)|\,dy.
	\end{aligned}
\end{equation}
Similarly, by \eqref{eq:1.10} we have
\begin{equation}
\label{eq:3.2}
\begin{split}
	 & \biggr|\biggr|\biggr|\,\int_{\{|y|<R^{\frac{1}{\theta}}\}}
 	|\nabla^j{\mathcal S}^m_K(\cdot,y,t)||\varphi(y)|\,dy\,\biggr|\biggr|\biggr|_{q,\ell}\\
	& \le C\biggr|\biggr|\biggr|\,\int_{\{|y|<R^{\frac{1}{\theta}}\}}
	\int_0^1|(\partial_t^m\nabla^{[K]+j+1}G_\theta)(\cdot-\tau y,t)||y|^{[K]+1}|\varphi(y)|\,d\tau\,dy\,\biggr|\biggr|\biggr|_{q,\ell}\\
	& \le C\int_{\{|y|<R^{\frac{1}{\theta}}\}}
	\int_0^1 |||(\partial_t^m\nabla^{[K]+j+1}G_\theta)(\cdot-\tau y,t)|||_{q,\ell}|y|^{[K]+1}|\varphi(y)|\,d\tau\,dy.
\end{split}
\end{equation}
On the other hand, it follows that 	
\begin{equation*}
\begin{split}
 |x|^\ell |(\partial_t^m\nabla^{[K]+j+1}G_\theta)(x-\tau y,t)|
 & =|z+\tau y|^\ell |(\partial_t^m\nabla^{[K]+j+1}G_\theta)(z,t)|\\
 & \le C(|z|^\ell+|y|^\ell)|(\partial_t^m\nabla^{[K]+j+1}G_\theta)(z,t)|
\end{split}
\end{equation*}
for $x$, $y\in{\mathbb R}^N$, $t>0$, and $\tau\in(0,1)$, where $z:=x-\tau y$. 
This together with \eqref{eq:2.2} implies that 
\begin{equation}
\label{eq:3.3}
\begin{split}
 & |||(\partial_t^m\nabla^{[K]+j+1}G_\theta)(\cdot-\tau y,t)|||_{q,\ell}\\
 & \le|||(\partial_t^m\nabla^{[K]+j+1}G_\theta)(t)|||_{q,\ell}
+|y|^\ell|||(\partial_t^m\nabla^{[K]+j+1}G_\theta)(t)|||_q\\
 & \le Ct^{-\frac{N}{\theta}\left(1-\frac{1}{q}\right)-m-\frac{[K]+j+1-\ell}{\theta}}
+Ct^{-\frac{N}{\theta}\left(1-\frac{1}{q}\right)-m-\frac{[K]+j+1}{\theta}}|y|^\ell
\end{split}
\end{equation}
for $y\in{\mathbb R}^N$, $t>0$, and $\tau\in(0,1)$.  
By \eqref{eq:3.2} and \eqref{eq:3.3} we obtain 
\begin{equation}
\label{eq:3.4}
	\begin{aligned}
	 & t^{\frac{N}{\theta}\left(1-\frac{1}{q}\right)+m+\frac{j}{\theta}}\biggr|\biggr|\biggr|\,\int_{\{|y|<R^{\frac{1}{\theta}}\}}
 	|\nabla^j{\mathcal S}^m_K(\cdot,y,t)||\varphi(y)|\,dy\,\biggr|\biggr|\biggr|_{q,\ell}\\
 	& \le C\int_{\{|y|<R^{\frac{1}{\theta}}\}} 
 	(t^{-\frac{[K]+1-\ell}{\theta}}+t^{-\frac{[K]+1}{\theta}} |y|^{\ell})|y|^{[K]+1} |\varphi(y)|\,dy\\
 	& \le C\left(t^{-\frac{[K]+1-\ell}{\theta}}R^{\frac{[K]+1-K}{\theta}}
 	+t^{-\frac{[K]+1}{\theta}}R^{\frac{[K]+\ell+1-K}{\theta}}\right)
 	\int_{\{|y|<R^{\frac{1}{\theta}}\}}|y|^K |\varphi(y)|\,dy.
	\end{aligned}
\end{equation}
Combining \eqref{eq:3.1} and \eqref{eq:3.4} and setting $R=t$, 
we obtain 
$$
	 t^{\frac{N}{\theta}\left(1-\frac{1}{q}\right)+m+\frac{j}{\theta}}\left|\left|\left|
	\int_{{\mathbb R}^N}\nabla^j{\mathcal S}^m_K(\cdot,y,t)\varphi(y)\,dy\right|\right|\right|_{q,\ell}
	\le Ct^{-\frac{K-\ell}{\theta}}|||\varphi|||_{1,K},\quad t>0,
$$
which implies assertion~(b). 
Similarly, setting $R=\varepsilon t$ with $0<\varepsilon\le 1$, 
we have 
\begin{align*}
 	& t^{\frac{N}{\theta}\left(1-\frac{1}{q}\right)+m+\frac{j}{\theta}}\left|\left|\left|
	\int_{{\mathbb R}^N}\nabla^j{\mathcal S}^m_K(\cdot,y,t)\varphi(y)\,dy\right|\right|\right|_{q,\ell}\\
 	& \le Ct^{-\frac{K-\ell}{\theta}}
	\biggr((\varepsilon^{-1})^{\frac{K-\ell}{\theta}}
	+\sum_{\ell<|\alpha|\le K}(\varepsilon^{-1})^{\frac{K-|\alpha|}{\theta}}\biggr)
	\int_{\{|y| \ge(\varepsilon t)^{\frac{1}{\theta}}\}}|y|^K|\varphi(y)|\,dy\\
 	&
	\qquad
	 +Ct^{-\frac{K-\ell}{\theta}}\left(\varepsilon^{\frac{[K]+1-K}{\theta}}
	+\varepsilon^{\frac{[K]+\ell+1-K}{\theta}}\right)|||\varphi|||_{1,K}.
\end{align*}
This together with $\varphi\in L^1_K$ implies that 
\begin{align*}
 	& \limsup_{t\to\infty}\,t^{\frac{N}{\theta}\left(1-\frac{1}{q}\right)+m+\frac{K+j-\ell}{\theta}}
 	\biggr|\biggr|\biggr|\,\int_{{\mathbb R}^N}\nabla^j{\mathcal S}^m_K(\cdot,y,t)\varphi(y)\,dy\,
	\biggr|\biggr|\biggr|_{q,\ell}\\
 	& \le C\left(\varepsilon^{\frac{[K]+1-K}{\theta}}+\varepsilon^{\frac{[K]+\ell+1-K}{\theta}}\right)|||\varphi|||_{1,K}. 
\end{align*}
Since $\varepsilon$ is arbitrary,  
we obtain assertion~(c). 
Thus Lemma~\ref{Lemma:3.2} follows. 
$\Box$
\vspace{5pt}

By Lemmata~\ref{Lemma:3.1} and \ref{Lemma:3.2} we have:
\begin{lemma}
\label{Lemma:3.3}
Let $f$ be a measurable function in $\mathbb R^N\times(0,\infty)$. 
Assume \eqref{eq:1.5} for some $K\ge 0$ and $1\le q\le\infty$. 
Let $0\le\ell\le K$, $m\in{\mathbb N}_0$, and $j=0,1$.  
Then there exists $C>0$ such that
\begin{equation}
\label{eq:3.5}
\begin{split}
	 & (t+1)^{\frac{N}{\theta}\left(1-\frac{1}{q}\right)+\frac{K+j-\ell}{\theta}}
	\biggr|\biggr|\biggr|\,\int_T^t\int_{{\mathbb R}^N}(s+1)^m\nabla^j{\mathcal S}^m_K(\cdot,y,t+1)f(y,s)\,dy\,ds\,
	\biggr|\biggr|\biggr|_{q,\ell}\\
	 & \le C\int_T^t E_{K,q}[f](s)\,ds
\end{split}
\end{equation}
for $t>T\ge 0$. Furthermore, 
$$
	\lim_{t\to\infty}
	t^{\frac{N}{\theta}\left(1-\frac{1}{q}\right)+\frac{K+j-\ell}{\theta}}
	\biggr|\biggr|\biggr|\,\int_0^T\int_{{\mathbb R}^N}(s+1)^m\nabla^j{\mathcal S}^m_K(\cdot,y,t+1)f(y,s)\,dy\,ds\,
	\biggr|\biggr|\biggr|_{q,\ell}=0
$$
for $T> 0$.
\end{lemma}
{\bf Proof.}
It follows from Lemma~\ref{Lemma:3.1} that 
\begin{equation}
\label{eq:3.6}
	|||f(s)|||_{1,K}\le CE_{K,q}[f](s),
	\quad s>0.
\end{equation}
This together with Lemma~\ref{Lemma:3.2}~(c) implies that, for any $T>0$, 
\begin{equation}
\label{eq:3.7}
	\lim_{t\to\infty}
	t^{\frac{N}{\theta}\left(1-\frac{1}{q}\right)+\frac{K+j-\ell}{\theta}}(s+1)^m
	\left|\left|\left|\,
	\int_{\mathbb R^N}\nabla^j{\mathcal S}^m_K(\cdot,y,t+1)f(y,s)\,dy\,\right|\right|\right|_{q,\ell}=0
\end{equation}
for $0<s<T$.
Furthermore, by \eqref{eq:3.6} with Lemma~\ref{Lemma:3.2}~(b) we see that
\begin{equation}
\label{eq:3.8}
	\begin{split}
 	& (t+1)^{\frac{N}{\theta}\left(1-\frac{1}{q}\right)+\frac{K+j-\ell}{\theta}}(s+1)^m\left|\left|\left|\,
	\int_{\mathbb R^N}\nabla^j{\mathcal S}^m_K(\cdot,y,t+1)f(y,s)\,dy\,\right|\right|\right|_{q,\ell}\\
 	& \le C(t+1)^{-m}(s+1)^m\|f(s)\|_{1,K}
	\le CE_{K,q}[f](s)
	\end{split}
\end{equation}
for $0<s<t$.
Inequality~\eqref{eq:3.8} implies \eqref{eq:3.5}.
Furthermore, by \eqref{eq:3.7} and \eqref{eq:3.8} we apply the Lebesgue dominated convergence theorem 
to obtain 
$$
	\lim_{t\to\infty}
	t^{\frac{N}{\theta}\left(1-\frac{1}{q}\right)+\frac{K+j-\ell}{\theta}}
	\biggr|\biggr|\biggr|\,
	\int_0^T\int_{{\mathbb R}^N}
	(s+1)^m \nabla^j{\mathcal S}_K^m(\cdot,y,t+1)f(y,s)\,dy\,ds
	\,\biggr|\biggr|\biggr|_{q,\ell}=0. 
$$
Thus Lemma~\ref{Lemma:3.3} follows.
$\Box$\vspace{5pt}

Next we prove the following lemma on the integral kernel ${\mathcal T}(x,y,t,s)$.
\begin{lemma}
\label{Lemma:3.4}
Let $1\le q\le\infty$, and $0\le \ell\le K$.
\begin{itemize}
  \item[{\rm (a)}] 
  Let $\varphi\in L^1_K$ with $K\ge 0$ and $j=0,1$. Then there exists $C_1>0$ such that 
	\begin{equation}
	\label{eq:3.9}
		t^{\frac{N}{q}\left(1-\frac{1}{q}\right)+\frac{j}{\theta}}(t+1)^{\frac{K-\ell}{\theta}}
		\biggr|\biggr|\biggr|\,\int_{{\mathbb R}^N}\nabla^j{\mathcal T}(\cdot,y,t+1,1)\varphi(y)\,dy\,
		\biggr|\biggr|\biggr|_{q,\ell}
		\le C_1\|\varphi\|_{L^1_K}
	\end{equation}
	for $t>0$.
	Furthermore,
	\begin{equation}
	\label{eq:3.10}
		\lim_{t\to\infty}t^{\frac{N}{\theta}\left(1-\frac{1}{q}\right)+\frac{K+j-\ell}{\theta}}
		\biggr|\biggr|\biggr|\,\int_{{\mathbb R}^N}\nabla^j{\mathcal T}(\cdot,y,t+1,1)\varphi(y)\,dy\, 
		\biggr|\biggr|\biggr|_{q,\ell}=0.
	\end{equation}
  \item[{\rm (b)}]  
  Let $f$ be a measurable function in $\mathbb R^N\times(0,\infty)$ 
  and satisfy \eqref{eq:1.5}. 
  Let  $j=0$ if $0<\theta\le 1$ and $j\in\{0,1\}$ if $1\le\theta<2$. 
  Then there exists $C_2>0$ such that 
		\begin{equation}
		\label{eq:3.11}
			\begin{aligned}
			&
			t^{\frac{N}{q}\left(1-\frac{1}{q}\right)}(t+1)^{\frac{K-\ell}{\theta}}
			\biggr|\biggr|\biggr|\,\int_T^t\int_{{\mathbb R}^N}\nabla^j {\mathcal T}(\cdot,y,t+1,s+1)f(y,s)\,dy\,ds\, 
			\biggr|\biggr|\biggr|_{q,\ell}\\
			& \le C_2\int_T^t (t-s)^{-\frac{j}{\theta}}E_{K,q}[f](s)\,ds
			\end{aligned}
		\end{equation}
		for $t> T\ge0$.
		Furthermore,
		\begin{equation}
		\label{eq:3.12}
			\lim_{t\to\infty}t^{\frac{N}{q}\left(1-\frac{1}{q}\right)+\frac{K-\ell}{\theta}}
			\biggr|\biggr|\biggr|\,\int_0^T\int_{{\mathbb R}^N}{\mathcal T}(\cdot,y,t+1,s+1)f(y,s)\,dy\,ds\, 
			\biggr|\biggr|\biggr|_{q,\ell}=0
		\end{equation}
		for $T> 0$.
\end{itemize}
\end{lemma}
{\bf Proof.} 
Let $0\le \ell\le K$ and $j=0,1$.
We find $\ell'>0$ such that
\begin{equation}
\label{eq:3.13}
	\ell\le \ell',
	\qquad
	\theta K_\theta<\ell'<\theta(K_\theta+1).
\end{equation}
Let $x$, $y\in{\mathbb R}^N$ and $t>0$. 
It follows that 
\begin{equation}
\label{eq:3.14}
	t^{-\frac{\ell}{\theta}}|x|^\ell
	\le t^{-\frac{\ell}{\theta}}(|x-y|^{\ell}+|y|^\ell)
	\le C\left(1+t^{-\frac{\ell'}{\theta}}|x-y|^{\ell'}+t^{-\frac{K}{\theta}}|y|^K\right).
\end{equation}
This together with \eqref{eq:1.10} implies that
\begin{equation}
\label{eq:3.15}
	\begin{aligned}
 	& 
	t^{-\frac{\ell}{\theta}}|x|^\ell |\nabla^j{\mathcal T}(x,y,t,s)|\\
	& 
	\le Ct^{-\frac{\ell'}{\theta}}s^{K_\theta+1}
	\int_0^1|x-y|^{\ell'}|(\partial_t^{K_\theta+1}\nabla^jG_\theta)(x-y,t-\tau s)|\,d\tau
	\\
	&\quad
	+C\bigg(1+t^{-\frac{K}{\theta}}|y|^K\bigg)
	\left[|(\nabla^jG_\theta)(x-y,t-s)|+\sum_{m=0}^{K_\theta} s^m|(\partial_t^m\nabla^j G_\theta)(x-y,t)|\right]
	\end{aligned}
\end{equation}
for $0<s<t$. 
Let $\psi\in L^{r_1}_K\cap L^{r_2}_K$ with $1\le r_1, r_2\le q$. Let $1\le r_i'\le \infty$ $(i=1,2)$ be such that 
$$
\frac{1}{q}=\frac{1}{r_i}+\frac{1}{r_i'}-1.
$$
Then we observe from the Young inequality, \eqref{eq:2.2} and \eqref{eq:3.15} that 
\begin{equation}
\label{eq:3.16}
	\begin{split}
 	& 
	t^{-\frac{\ell}{\theta}}
	\biggr|\biggr|\biggr|\int_{{\mathbb R}^N}\nabla^j{\mathcal T}(\cdot,y,t,s)\psi(y)\,dy\,
	\biggr|\biggr|\biggr|_{q,\ell}\\
 	& 
	\le Ct^{-\frac{\ell'}{\theta}}s^{K_\theta+1}
	\int_0^1|||\partial_t^{K_\theta+1}\nabla^jG_\theta(t-\tau s)|||_{r'_1,\ell'}\|\psi\|_{r_1}\,d\tau
	\\
	&\qquad
	+C\|\nabla^jG_\theta(t-s)\|_{r'_1}\|\varphi\|_{r_1}
	+Ct^{-\frac{K}{\theta}}\|\nabla^jG_\theta(t-s)\|_{r'_1}|||\varphi|||_{r_1,K}\\
 	& \qquad\qquad
	+C\sum_{m=0}^{K_\theta} s^m\bigg[\|\partial_t^m \nabla^jG_\theta(t)\|_{r'_2}\|\varphi\|_{r_2}
	+t^{-\frac{K}{\theta}}\|\partial_t^m\nabla^j G_\theta(t)\|_{r'_2}|||\varphi|||_{r_2,K}\bigg]\\
	& 
	\le Ct^{-\frac{\ell'}{\theta}}s^{K_\theta+1}\|\psi\|_{r_1}
	\int_0^1(t-\tau s)^{-\frac{N}{\theta}\left(\frac{1}{r_1}-\frac{1}{q}\right)-(K_\theta+1)+\frac{\ell'-j}{\theta}}\,d\tau
	\\
	&\qquad
	+C(t-s)^{-\frac{N}{\theta}\left(\frac{1}{r_1}-\frac{1}{q}\right)-\frac{j}{\theta}}
	(\|\psi\|_{r_1}+t^{-\frac{K}{\theta}}|||\psi|||_{r_1,K})\\
	& \qquad\qquad
	+Ct^{-\frac{N}{\theta}\left(\frac{1}{r_2}-\frac{1}{q}\right)-\frac{j}{\theta}}
	\bigg(1+\sum_{m=0}^{K_\theta} s^mt^{-m}\bigg)
	(\|\psi\|_{r_2}+t^{-\frac{K}{\theta}}|||\psi|||_{r_2,K})
	\end{split}
\end{equation}
for $t/2< s<t$.
On the other hand, it follows from \eqref{eq:3.13} that 
\begin{equation*}
\begin{split}
\int_0^1(t-\tau s)^{-\frac{N}{\theta}\left(\frac{1}{r_1}-\frac{1}{q}\right)-(K_\theta+1)+\frac{\ell'-j}{\theta}}\,d\tau
 & \le (t-s)^{-\frac{N}{\theta}\left(\frac{1}{r_1}-\frac{1}{q}\right)}
\int_0^1(t-\tau s)^{-(K_\theta+1)+\frac{\ell'-j}{\theta}}\,d\tau\\
 & \le C(t-s)^{-\frac{N}{\theta}\left(\frac{1}{r_1}-\frac{1}{q}\right)-\frac{j}{\theta}}s^{-1}t^{-K_\theta+\frac{\ell'}{\theta}}
\end{split}
\end{equation*}
for $0<s<t$. 
This together with \eqref{eq:3.16} implies that 
\begin{equation}
\label{eq:3.17}
	\begin{split}
	 & t^{-\frac{\ell}{\theta}}
	\biggr|\biggr|\biggr|\int_{{\mathbb R}^N}\nabla^j{\mathcal T}(\cdot,y,t,s)\psi(y)\,dy\,
	\biggr|\biggr|\biggr|_{q,\ell}\\
	&
	\le C(t-s)^{-\frac{N}{\theta}\left(\frac{1}{r_1}-\frac{1}{q}\right)-\frac{j}{\theta}}
	(\|\psi\|_{r_1}+t^{-\frac{K}{\theta}}|||\psi|||_{r_1,K})\\
	& \qquad\qquad
	+Ct^{-\frac{N}{\theta}\left(\frac{1}{r_2}-\frac{1}{q}\right)}(t-s)^{-\frac{j}{\theta}}
	(\|\psi\|_{r_2}+t^{-\frac{K}{\theta}}|||\psi|||_{r_2,K})
	\end{split}
\end{equation}
for $t/2< s<t$.
Similarly, 
by \eqref{eq:1.10} and \eqref{eq:3.14} 
we have 
\begin{equation}
\label{eq:3.18}
	\begin{aligned}
	t^{-\frac{\ell}{\theta}}|x|^\ell |\nabla^j{\mathcal T}(x,y,t,s)
	 & |\le Cs^{K_\theta+1}\left(1+t^{-\frac{\ell'}{\theta}}|x-y|^{\ell'}+t^{-\frac{K}{\theta}}|y|^K\right)\\
	 & \quad
	 \times\int_0^1|(\partial_t^{K_\theta+1}\nabla^jG_\theta)(x-y,t-\tau s)|\,d\tau
	\end{aligned}
\end{equation}
for $0<s<t$. 
It follows from the Young inequality, \eqref{eq:2.2}, and \eqref{eq:3.18} that
\begin{equation}
\label{eq:3.19}
	\begin{split}
 	& 
	t^{-\frac{\ell}{\theta}}
	\biggr|\biggr|\biggr|\int_{{\mathbb R}^N}\nabla^j{\mathcal T}(\cdot,y,t,s)\psi(y)\,dy\,
	\biggr|\biggr|\biggr|_{q,\ell}\\
 	& 
	\le Cs^{K_\theta+1}
	\int_0^1\|(\partial_t^{K_\theta+1}\nabla^j G_\theta)(t-\tau s)\|_{r'_1}\|\psi\|_{r_1}\,d\tau\\
	& \qquad
	+Cs^{K_\theta+1}t^{-\frac{\ell'}{\theta}}
	\int_0^1|||(\partial_t^{K_\theta+1}\nabla^j G_\theta)(t-\tau s)|||_{r'_1,\ell'}\|\psi\|_{r_1}\,d\tau
	\\
	& \qquad\qquad
	+Cs^{K_\theta+1}t^{-\frac{K}{\theta}}
	\int_0^1\|(\partial_t^{K_\theta+1}\nabla^j G_\theta)(t-\tau s)\|_{r'_1}|||\psi|||_{r_1,K}\,d\tau\\
	&
	\le Cs^{K_\theta+1}\|\psi\|_{r_1}
	\int_0^1(t-\tau s)^{-\frac{N}{\theta}\left(\frac{1}{r_1}-\frac{1}{q}\right)-(K_\theta+1)-\frac{j}{\theta}}\,d\tau\\
	& \qquad
	+Cs^{K_\theta+1}t^{-\frac{\ell'}{\theta}}\|\psi\|_{r_1}
	\int_0^1(t-\tau s)^{-\frac{N}{\theta}\left(\frac{1}{r_1}-\frac{1}{q}\right)-(K_\theta+1)+\frac{\ell'-j}{\theta}}\,
	d\tau\\
	& \qquad\qquad
	+Cs^{K_\theta+1}t^{-\frac{K}{\theta}}|||\psi|||_{r_1,K}
	\int_0^1(t-\tau s)^{-\frac{N}{\theta}\left(\frac{1}{r_1}-\frac{1}{q}\right)-(K_\theta+1)-\frac{j}{\theta}}\,d\tau\\
 	& \le Cs^{K_\theta+1}t^{-\frac{N}{\theta}\left(\frac{1}{r_1}-\frac{1}{q}\right)-(K_\theta+1)}
	(t-s)^{-\frac{j}{\theta}}
	\bigg\{\|\psi\|_{r_1}+t^{-\frac{K}{\theta}}|||\psi|||_{r_1,K}\bigg\}\\
	\end{split}
\end{equation}
for $0<s\le t/2$.

We prove assertion~(a). 
Since $(t+1)/2\le 1$ for $0<t\le 1$, 
by \eqref{eq:3.17} with $\psi=\varphi$ and $r_1=r_2=1$ we have 
\begin{equation}
\label{eq:3.20}
	\begin{split}
 	& 
	(t+1)^{-\frac{\ell}{\theta}}
	\biggr|\biggr|\biggr|\int_{{\mathbb R}^N}\nabla^j{\mathcal T}(\cdot,y,t+1,1)\varphi(y)\,dy\,
	\biggr|\biggr|\biggr|_{q,\ell}\\
	 & 
	 \le Ct^{-\frac{N}{\theta}\left(1-\frac{1}{q}\right)-\frac{j}{\theta}}
	 (\|\varphi\|_1+(t+1)^{-\frac{K}{\theta}}|||\varphi|||_{1,K})\\
	 &
	 \le Ct^{-\frac{N}{\theta}\left(1-\frac{1}{q}\right)-\frac{j}{\theta}}\|\varphi\|_{L^1_K}
	\end{split}
\end{equation}
for $0<t\le 1$.
On the other hand, since  $(t+1)/2>1$ for $t>1$
by \eqref{eq:3.19} with $\psi=\varphi$ and $r_1=1$ we see that 
\begin{equation*}
	\begin{split}
 	& (t+1)^{-\frac{\ell}{\theta}}
	\biggr|\biggr|\biggr|\int_{{\mathbb R}^N}\nabla^j{\mathcal T}(\cdot,y,t+1,1)\varphi(y)\,dy\,
	\biggr|\biggr|\biggr|_{q,\ell}\\
 	& \le C(t+1)^{-\frac{N}{\theta}\left(1-\frac{1}{q}\right)-(K_\theta+1)}t^{-\frac{j}{\theta}}
	\bigg\{\|\varphi\|_1+(t+1)^{-\frac{K}{\theta}}|||\varphi|||_{1,K}\bigg\}\\
	& \le Ct^{-\frac{N}{\theta}\left(1-\frac{1}{q}\right)-\frac{j}{\theta}}(t+1)^{-(K_\theta+1)}\|\varphi\|_{L^1_K}
	\end{split}
\end{equation*}
for $t>1$.
This together with
$\theta(K_\theta+1)>K$ implies \eqref{eq:3.10} and 
\begin{equation}
\label{eq:3.21}
	t^{\frac{N}{\theta}\left(1-\frac{1}{q}\right)+\frac{j}{\theta}}(t+1)^{\frac{K-\ell}{\theta}}
 	\biggr|\biggr|\biggr|\int_{{\mathbb R}^N}{\mathcal T}(\cdot,y,t+1,1)\varphi(y)\,dy\,\biggr|\biggr|\biggr|_{q,\ell}
	\le C\|\varphi\|_{L^1_K}
\end{equation}
for $t> 1$. 
Combining \eqref{eq:3.20} and \eqref{eq:3.21}, we obtain \eqref{eq:3.9} and \eqref{eq:3.10}. 
Thus assertion~(a) follows.

We prove assertion~(b).
Since $(t+1)/2<s+1$ for $t/2<s<t$,
by \eqref{eq:3.13} and \eqref{eq:3.17} with $\psi=f(s)$ and $(r_1,r_2)=(q,1)$ we have
$$
	\begin{aligned}
 	& 
	(t+1)^{-\frac{\ell}{\theta}}
	\biggr|\biggr|\biggr|\int_{{\mathbb R}^N}\nabla^j{\mathcal T}(\cdot,y,t+1,s+1)f(y,s)\,dy\,
	\biggr|\biggr|\biggr|_{q,\ell}\\
	 &
	 \le C(t-s)^{-\frac{j}{\theta}}
	(\|f(s)\|_q+(t+1)^{-\frac{K}{\theta}}|||f(s)|||_{q,K})\\
	& \qquad
	+Ct^{-\frac{N}{\theta}\left(1-\frac{1}{q}\right)}(t-s)^{-\frac{j}{\theta}}
	(\|f(s)\|_1+(t+1)^{-\frac{K}{\theta}}|||f(s)|||_{1,K})
	\end{aligned}
$$
for $t/2<s<t$.
This together with Lemma~\ref{Lemma:3.1} implies that
\begin{equation}
\label{eq:3.22}
\begin{split}
 	& 
	(t+1)^{-\frac{\ell}{\theta}}
	\biggr|\biggr|\biggr|\int_{{\mathbb R}^N}\nabla^j{\mathcal T}(\cdot,y,t+1,s+1)f(y,s)\,dy\,\biggr|\biggr|\biggr|_{q,\ell}\\
 	& 
	\le C\left(s^{-\frac{N}{\theta}\left(1-\frac{1}{q}\right)}+t^{-\frac{N}{\theta}\left(1-\frac{1}{q}\right)}\right)
	\left((s+1)^{-\frac{K}{\theta}}+(t+1)^{-\frac{K}{\theta}}\right)(t-s)^{-\frac{j}{\theta}}E_{K,q}[f](s)\\
	& 
	\le Ct^{-\frac{N}{\theta}\left(1-\frac{1}{q}\right)}(t+1)^{-\frac{K}{\theta}}(t-s)^{-\frac{j}{\theta}}E_{K,q}[f](s)
\end{split}
\end{equation}
for $t/2<s<t$.
On the other hand,
by \eqref{eq:3.13} and \eqref{eq:3.17} with $\psi=f(s)$ and $r_1=r_2=1$ we have
$$
	\begin{aligned}
 	& 
	(t+1)^{-\frac{\ell}{\theta}}
	\biggr|\biggr|\biggr|\int_{{\mathbb R}^N}\nabla^j{\mathcal T}(\cdot,y,t+1,s+1)f(y,s)\,dy\,
	\biggr|\biggr|\biggr|_{q,\ell}\\
	&
	 \le Ct^{-\frac{N}{\theta}\left(1-\frac{1}{q}\right)}(t-s)^{-\frac{j}{\theta}}
	(\|f(s)\|_1+(t+1)^{-\frac{K}{\theta}}|||f(s)|||_{1,K})
	\end{aligned}
$$
for $0<s\le t/2$ with $(t+1)/2< s+1$.
This together with Lemma~\ref{Lemma:3.1} implies that
\begin{equation}
\label{eq:3.23}
	\begin{aligned}
 	& 
	(t+1)^{-\frac{\ell}{\theta}}
	\biggr|\biggr|\biggr|\int_{{\mathbb R}^N}\nabla^j{\mathcal T}(\cdot,y,t+1,s+1)f(y,s)\,dy\,
	\biggr|\biggr|\biggr|_{q,\ell}\\
 	& 
	\le Ct^{-\frac{N}{\theta}\left(1-\frac{1}{q}\right)}
	\left((s+1)^{-\frac{K}{\theta}}+(t+1)^{-\frac{K}{\theta}}\right)
	(t-s)^{-\frac{j}{\theta}}E_{K,q}[f](s)
	\\
	&
	\le Ct^{-\frac{N}{\theta}\left(1-\frac{1}{q}\right)}(t+1)^{-\frac{K}{\theta}}
	(t-s)^{-\frac{j}{\theta}}E_{K,q}[f](s)
	\end{aligned}
\end{equation}
for $0<s\le t/2$ with $(t+1)/2< s+1$. 
Furthermore,
by \eqref{eq:3.19} with $\psi=f(s)$ and $r_1=1$ 
we have
$$
	\begin{aligned}
 	& 
	(t+1)^{-\frac{\ell}{\theta}}
	\biggr|\biggr|\biggr|\int_{{\mathbb R}^N}\nabla^j{\mathcal T}(\cdot,y,t+1,s+1)f(y,s)\,dy\,
	\biggr|\biggr|\biggr|_{q,\ell}\\
 	& 
	\le C(s+1)^{K_\theta+1}
	(t+1)^{-\frac{N}{\theta}\left(1-\frac{1}{q}\right)-(K_\theta+1)}(t-s)^{-\frac{j}{\theta}}
	\times
	\\
	&\hspace{5cm}
	\times
	\bigg\{\|f(s)\|_1+(t+1)^{-\frac{K}{\theta}}|||f(s)|||_{1,K}\bigg\}
	\end{aligned}
$$
for $0<s\le t/2$ with $(t+1)/2\ge s+1$.
This together with Lemma~\ref{Lemma:3.1} again implies that
\begin{equation}
\label{eq:3.24}
	\begin{aligned}
 	& 
	(t+1)^{-\frac{\ell}{\theta}}
	\biggr|\biggr|\biggr|\int_{{\mathbb R}^N}\nabla^j{\mathcal T}(\cdot,y,t+1,s+1)f(y,s)\,dy\,
	\biggr|\biggr|\biggr|_{q,\ell}\\
 	& 
	\le Ct^{-\frac{N}{\theta}\left(1-\frac{1}{q}\right)}(s+1)^{K_\theta+1}
	(t+1)^{-(K_\theta+1)}\left((s+1)^{-\frac{K}{\theta}}+(t+1)^{-\frac{K}{\theta}}\right)\times
	\\
	&\hspace{5cm}
	\times
	(t-s)^{-\frac{j}{\theta}}E_{K,q}[f](s)
	\\
	&
	\le Ct^{-\frac{N}{\theta}\left(1-\frac{1}{q}\right)}(s+1)^{K_\theta+1-\frac{K}{\theta}}
	(t+1)^{-(K_\theta+1)}(t-s)^{-\frac{j}{\theta}}E_{K,q}[f](s)
	\end{aligned}
\end{equation}
for $0<s\le t/2$ with $(t+1)/2\ge s+1$.
Then, by \eqref{eq:3.24}, for any $T>0$, we observe from $\theta(K_\theta+1)>K$ that
\begin{equation}
\label{eq:3.25}
	\lim_{t\to\infty}t^{\frac{N}{\theta}\left(1-\frac{1}{q}\right)+\frac{j}{\theta}}(t+1)^{\frac{K-\ell}{\theta}}
 	\biggr|\biggr|\biggr|\int_{{\mathbb R}^N}\nabla^j{\mathcal T}(\cdot,y,t+1,s+1)f(y,s)\,dy\,\biggr|\biggr|\biggr|_{q,\ell}=0
\end{equation}
for $0<s<T$.
Furthermore, by \eqref{eq:3.22}, \eqref{eq:3.23}, and \eqref{eq:3.24} we see that 
\begin{equation}
\label{eq:3.26}
	\begin{aligned}
	&
	t^{\frac{N}{\theta}\left(1-\frac{1}{q}\right)}(t+1)^{\frac{K-\ell}{\theta}}
 	\biggr|\biggr|\biggr|\int_{{\mathbb R}^N}\nabla^j{\mathcal T}(\cdot,y,t+1,s+1)f(y,s)\,dy\,\biggr|\biggr|\biggr|_{q,\ell}
	\\
	&
	\le C(t-s)^{-\frac{j}{\theta}}E_{K,q}[f](s)
	\end{aligned}
\end{equation}
for $0<s< t$. 
This implies \eqref{eq:3.11}. 
Furthermore, by \eqref{eq:3.25} and \eqref{eq:3.26} 
we apply the Lebesgue dominated convergence theorem 
to obtain \eqref{eq:3.12}. 
Thus assertion~(b) follows. 
The proof of Lemma~\ref{Lemma:3.4} is complete.
$\Box$\vspace{5pt}

Now we are ready to complete the proof of Theorem~\ref{Theorem:1.1}.
\vspace{5pt}
\newline
{\bf Proof of Theorem~\ref{Theorem:1.1}.}
Let $u$ and $w$ be as in Theorem~\ref{Theorem:1.1}. 
Then, by \eqref{eq:1.11} and \eqref{eq:1.12} we have

$$
	\begin{aligned}
 	& u(x,t)-w(x,t)\\
 	& =\int_{\mathbb R^N}{\mathcal T}(x,y,t+1,1)\varphi(y)\, dy
	+\int_0^t\int_{{\mathbb R}^N} {\mathcal T}(x,y,t+1,s+1)f(y,s)\,dy\,ds\\
	&
	\qquad
	+\sum_{m=0}^{K_\theta}\frac{(-1)^m}{m!}\int_{{\mathbb R}^N}{\mathcal S}^m_K(x,y,t+1)\varphi(y)\,dy\\
 	&
	\qquad\qquad
	+\sum_{m=0}^{K_\theta}\frac{(-1)^m}{m!}\int_0^t\int_{{\mathbb R}^N} (s+1)^m{\mathcal S}_K^m(x,y,t+1)f(y,s)\,dy\,ds.
	\end{aligned}
$$
We apply Lemmata~\ref{Lemma:3.2}, \ref{Lemma:3.3}, and \ref{Lemma:3.4} 
to obtain \eqref{eq:1.6} and \eqref{eq:1.8}. 
Then we easily see that \eqref{eq:1.9} holds. 
Thus Theorem~\ref{Theorem:1.1} follows. 
$\Box$\vspace{3pt}
\newline
{\bf Proof of Corollary~\ref{Corollary:1.1}.}
By property (G)-(i) we see that 
$g_{\alpha,m}(0)\in L^q_K$ is equivalent to $g_{\alpha,m}(t)\in L^q_K$ for $t\ge0$. 
Then Corollary~\ref{Corollary:1.1} follows from Theorem~\ref{Theorem:1.1}.
$\Box$
\begin{remark}
\label{Remark:3.1}
	{\rm(i)}
	The arguments of \cites{IIK,IK,IKM} are in the frameworks of $L^q$ and $L^1_K$. 
	On the other hand, the arguments in the proof of Theorem~{\rm\ref{Theorem:1.1}} 
	are in the framework of $L^q_K$. 
	This improvement enables us to 
	obtain HOAE of solutions to the Cauchy problem 
	for nonlinear fractional diffusion equations such as \eqref{eq:1.2}. 
	See Section~{\rm 5}.
	\vspace{3pt}
	\newline
	{\rm(ii)}
	Let $0\le K<\theta$ and $\ell=0$. 
	By similar arguments to those in the proof of Theorem~$\ref{Theorem:1.1}$ 
	we see that Theorem~{\rm\ref{Theorem:1.1}} holds with $E_{K,q}[f]$ replaced by  
	$$
	E_{K,q}'[f](t)
	 := (t+1)^{\frac{K}{\theta}}
	\left[t^{\frac{N}{\theta}\left(1-\frac{1}{q}\right)}\|f(t)\|_q+\|f(t)\|_1\right]+|||f(t)|||_{1,K}.
	$$
	See also \cite{IKK02}*{Theorem~1.2}.
\end{remark}
%

\section{Fractional convection-diffusion equation}
In this section 
we consider the Cauchy problem for a convection type inhomogeneous fractional diffusion equation
\begin{equation}
\label{eq:4.1}
	\partial_t u+(-\Delta)^{\frac{\theta}{2}}u={\rm div} f(x,t)\quad\mbox{in}\quad{\mathbb R}^N\times(0,\infty),
	\quad
	u(x,0)=\varphi(x)\quad\mbox{in}\quad{\mathbb R}^N,
\end{equation}
where $1<\theta<2$, $\varphi\in L^1_K$ with $K\ge 0$, and $f=(f_1,\dots,f_N)$ is a vector-valued function 
in ${\mathbb R}^N\times(0,\infty)$.
Similarly to Theorem~\ref{Theorem:1.1}, we have:
\begin{theorem}
\label{Theorem:4.1}
	Let $N\ge1$, $1<\theta<2$, $K\ge 0$, and $1\le q\le\infty$. 
	Let $f=(f_1,\dots,f_N)$ be a vector-valued measurable function in $\mathbb R^N\times(0,\infty)$ satisfying \eqref{eq:1.5}.
	Let $u\in C({\mathbb R}^N\times(0,\infty))$ be a solution to problem~\eqref{eq:4.1}, that is, 
	$u$ satisfies
	$$
		u(x,t)
		=\int_{{\mathbb R}^N}G_\theta(x-y,t)\varphi(y)\,dy
		+\int_0^t\int_{{\mathbb R}^N}\nabla G_\theta(x-y,t-s)\cdot f(y,s)\,dy\,ds
	$$
	for $(x,t)\in{\mathbb R}^N\times(0,\infty)$, where $\varphi\in L^1_K$. 
	Let $0\le \ell\le K$. 
	Then
	\begin{equation}
	\label{eq:4.2}
 		\sup_{0<t<\tau}\,t^{\frac{N}{\theta}\left(1-\frac{1}{q}\right)}|||u(t)-z(t)|||_{q,\ell}<\infty
		\quad\mbox{for}\quad \tau>0,
	\end{equation}
	where 
	\begin{equation*}
	\begin{split}
		z(x,t) & :=
		\sum_{m=0}^{K_\theta}\sum_{|\alpha|\le K}
		M_\alpha(\varphi)\,g_{\alpha,m}(x,t)\\
		 & \qquad
		 +\sum_{m=0}^{K_\theta}\sum_{|\alpha|\le K}
		\sum_{j=1}^N\left(\int_0^t(s+1)^m M_\alpha(f_j(s))\,ds\right)\partial_{x_j}g_{\alpha,m}(x,t).
	\end{split}
	\end{equation*}
	Furthermore, there exists $C>0$ such that, 
	for any $\varepsilon>0$ and $T>0$, 
	\begin{equation}
	\label{eq:4.3}
 	 	t^{\frac{N}{\theta}\left(1-\frac{1}{q}\right)-\frac{\ell}{\theta}}|||u(t)-z(t)|||_{q,\ell}\\
  	 	\le\varepsilon t^{-\frac{K}{\theta}}
		+Ct^{-\frac{K}{\theta}}\int_T^t(t-s)^{-\frac{1}{\theta}}E_{K,q}[f](s)\,ds
	\end{equation}
	holds for large enough $t>0$.
\end{theorem}
{\bf Proof of Theorem~\ref{Theorem:4.1}.}
Let $u$ and $z$ be as in Theorem~\ref{Theorem:4.1}. 
Then, similarly to \eqref{eq:1.11} and \eqref{eq:1.12}, we have
$$
	\begin{aligned}
 	& 
	u(x,t)-z(x,t)\\
 	& 
	=\int_{\mathbb R^N}{\mathcal T}(x,y,t+1,1)\varphi(y)\, dy
	+\int_0^t\int_{{\mathbb R}^N} \nabla {\mathcal T}(x,y,t+1,s+1)\cdot f(y,s)\,dy\,ds\\
	& \quad
	+\sum_{m=0}^{K_\theta}\frac{(-1)^m}{m!}
	\int_{{\mathbb R}^N}{\mathcal S}^m_K(x,y,t+1)\varphi(y)\,dy\\
 	& \quad
	+\sum_{m=0}^{K_\theta}\frac{(-1)^m}{m!}\int_0^t\int_{{\mathbb R}^N}(s+1)^m\nabla {\mathcal S}_K^m(x,y,t+1)\cdot f(y,s)\,dy\,ds.
	\end{aligned}
$$
Similarly to the proof of Theorem~\ref{Theorem:1.1}, 
we apply Lemmata~\ref{Lemma:3.2}, \ref{Lemma:3.3}, and \ref{Lemma:3.4} 
to obtain \eqref{eq:4.2} and \eqref{eq:4.3}. 
Thus Theorem~\ref{Theorem:4.1} follows. 
$\Box$
\section{Nonlinear fractional diffusion equation}
Let $N\ge 1$, $0<\theta<2$, and $F\in C({\mathbb R}^N\times[0,\infty)\times{\mathbb R})$. 
Consider the Cauchy problem 
for a nonlinear fractional diffusion equation
\begin{equation}
\tag{P}
	\partial_t u+(-\Delta)^{\frac{\theta}{2}}u=F(x,t,u)\quad\mbox{in}\quad{\mathbb R}^N\times(0,\infty),
	\quad
	u(x,0)=\varphi(x)\quad\mbox{in}\quad{\mathbb R}^N,
\end{equation}
where $\varphi\in L^1_K\cap L^\infty$ for some $K\ge 0$
under the following condition~(F):
\begin{itemize}
\item[(F)]
	there exists $p>1+\theta/N$ such that 
  	$$
  		|F(x,t,v)-F(x,t,w)|\le C(|v|+|w|)^{p-1}|v-w|
  	$$
  	for $(x,t,v,w)\in{\mathbb R}^N\times[0,\infty)\times{\mathbb R}^2$.
\end{itemize}
Let $u\in C({\mathbb R}^N\times(0,\infty))$ be a solution to problem~(P) 
that is, $u$ satisfies 
$$
	u(x,t)=\left[e^{-t(-\Delta)^{\theta/2}}\varphi\right](x)
	+\int_0^t \left[e^{-(t-s)(-\Delta)^{\theta/2}}F(\cdot,s,u(\cdot,s))\right](x)\,ds
$$
for $(x,t)\in{\mathbb R}^N\times(0,\infty)$. 
In this section, under condition~\eqref{eq:1.13}, 
we obtain HOAE of the solution~$u$. 
Theorem~\ref{Theorem:5.1} is an application of Theorem~\ref{Theorem:1.1}.
\begin{theorem}
\label{Theorem:5.1}
	Let $N\ge 1$, $0<\theta<2$, and $\varphi\in L^1_K$ with $K\ge 0$. 
	Assume condition~{\rm (F)}. 
	Let $u\in C({\mathbb R}^N\times(0,\infty))$ be a solution to problem~{\rm (P)}. 
	Set 
	$$
	F(x,t):=F(x,t,u(x,t)),\quad (x,t)\in{\mathbb R}^N\times(0,\infty).
	$$ 
	\begin{itemize}
	\item[\rm (a)]
		Assume that $\varphi\in L^\infty_k$ with $k=\min\{N+\theta,K\}$. 
		Let $u$ satisfy \eqref{eq:1.13}.
  		Then
  		\begin{equation}
		\label{eq:5.1}
  		\sup_{t>0}\,(t+1)^{\frac{N}{q}\left(1-\frac{1}{q}\right)-\frac{\ell}{\theta}}|||u(t)|||_{q,\ell}<\infty
  		\end{equation}
  		for $1\le q\le\infty$ and $0\le\ell\le K$ with $\ell<\theta+N(1-1/q)$.
	\item[\rm (b)]
		Let $u$ satisfy \eqref{eq:5.1}. 
		If $p(N+\theta)>K+N$, 
		then
		\begin{equation}
		\label{eq:5.2}
			E_{K,q}[F](t)\le C(t+1)^{\frac{K}{\theta}-A_p},\quad t>0,
		\end{equation}
		for $1\le q\le\infty$, where $A_p:=N(p-1)/\theta>1$.
	\item[\rm(c)]
		Assume that \eqref{eq:5.2} holds. 
		Set
		\begin{equation}
		\label{eq:5.3}
			U_0(x,t):=\sum_{m=0}^{K_\theta}\sum_{|\alpha|\le K}
			\left(M_\alpha(\varphi)+\int_0^t (s+1)^m M_\alpha(F(s))\,ds\right)\,g_{\alpha,m}(x,t).
		\end{equation}
		Then
		\begin{equation}
		\label{eq:5.4}
			\sup_{t>0}\,t^{\frac{N}{\theta}\left(1-\frac{1}{q}\right)}(t+1)^{-\frac{\ell}{\theta}}
			|||u(t)-U_0(t)|||_{q,\ell}<\infty
		\end{equation}
		and
		\begin{equation}
		\label{eq:5.5}
			t^{\frac{N}{\theta}\left(1-\frac{1}{q}\right)-\frac{\ell}{\theta}}|||u(t)-U_0(t)|||_{q,\ell}
			=\left\{
			\begin{array}{ll}
			o(t^{-\frac{K}{\theta}})+O(t^{-A_p+1})  & \mbox{if}\quad A_p-1\not=K/\theta,\vspace{5pt}\\
			O(t^{-\frac{K}{\theta}}\log t)  & \mbox{if}\quad A_p-1=K/\theta,\vspace{3pt}
			\end{array}
			\right.
		\end{equation}
		as $t\to\infty$, for $1\le q\le\infty$ and $0\le \ell \le K$.  
	\end{itemize}
\end{theorem}
\begin{remark}
\label{Remark:5.1}
	Let $N\ge 1$, $0<\theta<2$, and $\varphi\in L^\infty$.  
	Assume condition~{\rm (F)}.
	%
	\begin{itemize}
	\item[\rm (i)] 
		There exists $\delta>0$ such that, if $\|\varphi\|_{L^{N(p-1)/\theta}}<\delta$, 
		then problem~{\rm (P)} possesses a solution~$u\in C({\mathbb R}^N\times(0,\infty))$ 
		satisfying \eqref{eq:1.13}. 
		See \cites{IKK01, IKO, W}. 
	%
	\item[\rm (ii)] 
		Let $F(x,t,u):=\lambda|u|^{p-1}u$ with $\lambda\le 0$. 
		Then the comparison principle implies that 
		$$
			|u(x,t)|\le \left[e^{-t(-\Delta)^{\theta/2}}|\varphi|\right](x),\quad(x,t)\in{\mathbb R}^N\times(0,\infty).
		$$
		This together with Lemma~{\rm\ref{Lemma:2.1}} implies \eqref{eq:1.13}.
	\end{itemize}
\end{remark}

We prepare the following lemma for the proof of Theorem~\ref{Theorem:5.1}.
\begin{lemma}
\label{Lemma:5.1}
	Assume condition~{\rm (F)}. Let $K\ge 0$. 
	Let $v_1$ and $v_2$ be measurable functions in ${\mathbb R}^N\times(0,\infty)$ 
	and $h$ in $(0,\infty)$ such that
	\begin{equation}
	\label{eq:5.6}
	\begin{split}
 		& (t+1)^{\frac{N}{\theta}\left(1-\frac{1}{q}\right)-\frac{\ell}{\theta}}|||v_i(t)|||_{q,\ell}<\infty,\quad i=1,2,\\
 		& (t+1)^{\frac{N}{\theta}\left(1-\frac{1}{q}\right)-\frac{\ell}{\theta}}|||v_1(t)-v_2(t)|||_{q,\ell}\le h(t),
	\end{split}
	\end{equation}
	for $t>0$, $1\le q\le\infty$, and $0\le\ell\le K$ with $\ell<\theta+N(1-1/q)$. 
	Assume that $p(N+\theta)>K+N$. Then 
	$$
		E_{K,q}[F(v_1)-F(v_2)](t)\le C(t+1)^{-A_p+\frac{K}{\theta}}h(t),\quad t>0. 
	$$
\end{lemma}
{\bf Proof.}
Let $0\le\ell\le K$ and $1\le q\le\infty$. 
Since $p(N+\theta)>K+N$, we find $\ell_1$, $\ell_2\ge 0$ such that  
$$
	0\le \ell_1<\theta+N,\quad 0\le \ell_2<\theta+N\left(1-\frac{1}{q}\right),
	\quad \ell=(p-1)\ell_1+\ell_2.
$$
Then, by condition~(F) and \eqref{eq:5.6} we see that
$$
	\begin{aligned}
 	& (t+1)^{\frac{N}{\theta}\left(1-\frac{1}{q}\right)-\frac{\ell}{\theta}}|||F(v_1(s))-F(v_2(s))|||_{q,\ell}\\
 	& \le C(t+1)^{\frac{N}{\theta}\left(1-\frac{1}{q}\right)-\frac{\ell}{\theta}}
 	(|||v_1(t)|||_{\infty,\ell_1}^{p-1}+|||v_2(t)|||_{\infty,\ell_1}^{p-1})|||v_1(t)-v_2(t)|||_{q,\ell_2}\\
 	& \le C(t+1)^{\frac{N}{\theta}\left(1-\frac{1}{q}\right)-\frac{\ell}{\theta}-\frac{N(p-1)}{\theta}
	+\frac{(p-1)\ell_1}{\theta}}|||v_1(t)-v_2(t)|||_{q,\ell_2}\\
 	& \le C(t+1)^{-\frac{\ell}{\theta}-\frac{N(p-1)}{\theta}+\frac{(p-1)\ell_1}{\theta}+\frac{\ell_2}{\theta}}h(t)
 	=C(t+1)^{-A_p}h(t),\quad t>0. 
	\end{aligned}
$$
Thus Lemma~\ref{Lemma:5.1} follows.
$\Box$\vspace{5pt}

\noindent
{\bf Proof of Theorem~\ref{Theorem:5.1}.} 
We prove assertion~(a). 
Since $A_p=N(p-1)/\theta>1$, 
the comparison principle together with condition~(F) and \eqref{eq:1.13} implies that 
$$
	|u(x,t)|\le \exp\left(C\int_0^t(s+1)^{-A_p}\,ds\right)\left[e^{-t(-\Delta)^{\theta/2}}|\varphi|\right](x)
	\le C\left[e^{-t(-\Delta)^{\theta/2}}|\varphi|\right](x)
$$
for $(x,t)\in{\mathbb R}^N\times(0,\infty)$. 
This together with Lemma~\ref{Lemma:2.1} implies assertion~(a). 
Furthermore, assertion~(b) follows from Lemma~\ref{Lemma:5.1} with $v_1=u$ and $v_2=0$. 
On the other hand,
by Theorem~\ref{Theorem:1.1} with \eqref{eq:5.2} we obtain \eqref{eq:5.4}. 
Furthermore, 
for any $\varepsilon>0$ and $T>0$, we have
$$
	t^{\frac{N}{\theta}\left(1-\frac{1}{q}\right)-\frac{\ell}{\theta}}|||u(t)-U_0(t)|||_{q,\ell}
	\le\varepsilon t^{-\frac{K}{\theta}}
	+Ct^{-\frac{K}{\theta}}\int_T^t(s+1)^{\frac{K}{\theta}-A_p}\,ds
$$
for large enough $t>0$.
This implies \eqref{eq:5.5}. Thus assertion~(c) follows. 
The proof of Theorem~\ref{Theorem:5.1} is complete. 
$\Box$\vspace{5pt}

As a corollary of Theorem~\ref{Theorem:5.1}, we have: 
\begin{Corollary}
\label{Corollary:5.1}
	Let $N\ge 1$, $0<\theta<2$, and $\varphi\in L_K^1$ with $K\ge 0$. 
	Assume condition~{\rm (F)} and
	\begin{equation}
	\label{eq:5.7}
		p>1+\frac{2K+\theta}{N}.
	\end{equation}
	Let $u\in C({\mathbb R}^N\times(0,\infty))$ be a solution to problem~{\rm (P)} 
	and satisfy \eqref{eq:5.1}. 
	Then there exists a set $\{M_{\alpha,m}\}\subset{\mathbb R}$, 
	where $m\in\{0,\dots,K_\theta\}$ and $\alpha\in{\mathbb M}$ with $|\alpha|\le K$, such that
	\begin{equation}
	\label{eq:5.8}
		t^{\frac{N}{\theta}\left(1-\frac{1}{q}\right)-\frac{\ell}{\theta}}|||u(t)-U_*(t)|||_{q,\ell}
		=o\left(t^{-\frac{K}{\theta}}\right)
		\quad\mbox{as}\quad t\to\infty
	\end{equation}
	for $1\le q\le\infty$ and $0\le \ell\le K$, 
	where
	\begin{equation}
	\label{eq:5.9}
		U_*(x,t):=\sum_{m=0}^{K_\theta}\sum_{|\alpha|\le K}M_{\alpha,m}\,g_{\alpha,m}(x,t). 
	\end{equation}
\end{Corollary}
{\bf Proof.}
It follows from \eqref{eq:5.7} that 
\begin{equation}
\label{eq:5.10}
	p(N+\theta)>K+N,\qquad
	A_p-1=\frac{N}{\theta}(p-1)-1>\frac{2K}{\theta}.
\end{equation}
By Theorem~\ref{Theorem:5.1} we have
\begin{equation}
\label{eq:5.11}
	t^{\frac{N}{\theta}\left(1-\frac{1}{q}\right)-\frac{\ell}{\theta}}
	|||u(t)-U_0(t)|||_{q,\ell}=o\left(t^{-\frac{K}{\theta}}\right)
	\quad\mbox{as}\quad t\to\infty
\end{equation}
for $1\le q\le\infty$ and $1\le \ell\le K$. 
Here $U_0$ is as in Theorem~\ref{Theorem:5.1}. 

Let $m\in\{0,\dots,K_\theta\}$ and $\alpha\in{\mathbb M}$ with $|\alpha|\le K$. 
Assertion~(b) of Theorem~\ref{Theorem:5.1} implies that
\begin{equation}
\label{eq:5.12}
	|M_\alpha(F(t))|\le C(t+1)^{-A_p+\frac{|\alpha|}{\theta}},\quad t>0.
\end{equation}
Then, by \eqref{eq:5.10} we find $M_{\alpha,m}\in {\mathbb R}$ such that
$$
	M_{\alpha,0}=M_\alpha(\varphi)+\int_0^\infty M_{\alpha}(F(s))\,ds,
	\quad
	M_{\alpha,m}=\int_0^\infty (s+1)^mM_{\alpha}(F(s))\,ds\quad (m\ge1).
$$
Furthermore, by \eqref{eq:2.2}, \eqref{eq:5.3}, \eqref{eq:5.9}, \eqref{eq:5.10}, and \eqref{eq:5.12} 
we have
$$
	\begin{aligned}
 	& t^{\frac{N}{\theta}\left(1-\frac{1}{q}\right)-\frac{\ell}{\theta}}|||U_0(t)-U_*(t)||_{q,\ell}\\
 	& \le t^{\frac{N}{\theta}\left(1-\frac{1}{q}\right)-\frac{\ell}{\theta}}
 	\sum_{m=0}^{K_\theta}\sum_{|\alpha|\le K}
	\left(\int_t^\infty (s+1)^m |M_\alpha(F(s))|\,ds\right)\,|||g_{\alpha,m}(t)|||_{q,\ell}\\
  	& \le C\sum_{m=0}^{K_\theta}\sum_{|\alpha|\le K}(t+1)^{-m-\frac{|\alpha|}{\theta}}
	\int_t^\infty (s+1)^m(s+1)^{-A_p+\frac{|\alpha|}{\theta}}\,ds
  	\le Ct^{-A_p+1},\quad t\ge 1.
	\end{aligned}
$$
This together with \eqref{eq:5.10} and \eqref{eq:5.11} implies \eqref{eq:5.8}.  
Thus Corollary~\ref{Corollary:5.1} follows.
$\Box$\vspace{5pt}

Combining Theorems~\ref{Theorem:1.1} and \ref{Theorem:5.1}, 
we obtain a refined asymptotic expansion of the solution to problem~(P).
\begin{theorem}
\label{Theorem:5.2}
	Assume the same conditions as in Theorem~$\ref{Theorem:5.1}$. 
	Let $u$ satisfy \eqref{eq:5.1}. 
	For $n=1,2,\dots$, 
	define a function $U_n=U_n(x,t)$ in ${\mathbb R}^N\times(0,\infty)$ inductively by
	$$
		\begin{aligned}
		U_n(x,t) & :=U_0(x,t)+\int_0^t \left[e^{-(t-s)(-\Delta)^{\theta/2}}F_{n-1}(s)\right](x)\,ds\\
 		& \quad-\sum_{m=0}^{K_\theta}\sum_{|\alpha|\le K}
		\left(\int_0^t (s+1)^m M_\alpha(F_{n-1}(s))\,ds\right)\,g_{\alpha,m}(x,t),
		\end{aligned}
	$$
	where $U_0$ is as in Theorem~{\rm\ref{Theorem:5.1}} and $F_{n-1}(x,t):=F(x,t,U_{n-1}(x,t))$.
	Then
	\begin{equation}
	\label{eq:5.13}
		\sup_{t>0}\,t^{\frac{N}{\theta}\left(1-\frac{1}{q}\right)-\frac{\ell}{\theta}}
		|||u(t)-U_n(t)|||_{q,\ell}<\infty
	\end{equation}
	and
	\begin{equation}
	\label{eq:5.14}
		\begin{aligned}
 		& t^{\frac{N}{\theta}\left(1-\frac{1}{q}\right)-\frac{\ell}{\theta}}|||u(t)-U_n(t)|||_{q,\ell}\\
 		&\qquad
		=\left\{
		\begin{array}{ll}
		o(t^{-\frac{K}{\theta}})+O(t^{-(n+1)(A_p-1)})  & \mbox{if}\quad (n+1)(A_p-1)\not=K/\theta,\vspace{5pt}\\
		O(t^{-\frac{K}{\theta}}\log t)  & \mbox{if}\quad (n+1)(A_p-1)=K/\theta,\\
		\end{array}
		\right.
		\end{aligned}
	\end{equation}
	as $t\to\infty$, for $1\le q\le\infty$ and $0\le \ell\le K$.
\end{theorem}
{\bf Proof.}
Let $K\ge 0$. 
By Theorem~\ref{Theorem:5.1} we have \eqref{eq:5.13} and \eqref{eq:5.14} with $n=0$.
Assume that \eqref{eq:5.13} and \eqref{eq:5.14} hold for some $n=k\in\{0,1,\dots\}$. 
Then, by \eqref{eq:5.1} and \eqref{eq:5.13} with $n=k$  we have
\begin{equation}
\label{eq:5.15}
	\begin{aligned}
 	& \sup_{t>0}\,t^{\frac{N}{\theta}\left(1-\frac{1}{q}\right)-\frac{\ell}{\theta}}|||U_k(t)|||_{q,\ell}\\
 	& \le \sup_{t>0}\,t^{\frac{N}{\theta}\left(1-\frac{1}{q}\right)-\frac{\ell}{\theta}}|||u(t)-U_k(t)|||_{q,\ell}
	+\sup_{t>0}\,t^{\frac{N}{\theta}\left(1-\frac{1}{q}\right)-\frac{\ell}{\theta}}|||u(t)|||_{q,\ell}<\infty
	\end{aligned}
\end{equation}
for $1\le q\le\infty$ and $0\le \ell\le K$ with $\ell<\theta+N(1-1/q)$. 
On the other hand, it follows that
\begin{equation}
\label{eq:5.16}
	\begin{aligned}
 	& u(x,t)-\int_0^t \left[e^{-(t-s)(-\Delta)^{\theta/2}}F_k(s)\right](x)\,ds\\
 	& =\left[e^{-t(-\Delta)^{\theta/2}}\varphi\right](x)
	+\int_0^te^{-(t-s)(-\Delta)^{\theta/2}}[F(s)-F_k(s)]\,ds
	\end{aligned}
\end{equation}
for $(x,t)\in{\mathbb R}^N\times(0,\infty)$. 
By \eqref{eq:5.2} and \eqref{eq:5.15} 
we apply Lemma~\ref{Lemma:5.1} to obtain 
$$
	E_{K,q}[F-F_k]\in L^\infty(0,\tau)\quad\mbox{for}\quad\tau>0
$$
and
$$
	\begin{aligned}
	E_{K,q}[F-F_k](t) 
 	& =\left\{
	\begin{array}{ll}
	o\left(t^{-A_p}\right) & \mbox{if}\quad (k+1)(A_p-1)<K/\theta,\vspace{5pt}\\
	O\left(t^{-A_p}\log t\right) & \mbox{if}\quad (k+1)(A_p-1)=K/\theta,\vspace{5pt}\\
	O\left(t^{-A_p+\frac{K}{\theta}-(k+1)(A_p-1)}\right) & \mbox{if}\quad (k+1)(A_p-1)>K/\theta,
	\end{array}
	\right.
	\end{aligned}
$$
as $t\to\infty$. 
Then we apply Theorem~\ref{Theorem:1.1} to \eqref{eq:5.16}, 
namely $f(x,t)=F(x,t)-F_k(x,t)$,
and obtain \eqref{eq:5.13} and \eqref{eq:5.14} with $n=k+1$. 
Therefore, by induction we obtain \eqref{eq:5.13} and \eqref{eq:5.14} for $n=0,1,2,\dots$.  
Thus Theorem~\ref{Theorem:5.2} follows.
$\Box$\vspace{5pt}

Similarly to the proof of Theorem~\ref{Theorem:5.2} for the case $n=1$,
we prove Theorem~\ref{Theorem:1.2}.
\vspace{5pt}
\newline
{\bf Proof of Theorem~\ref{Theorem:1.2}.}
Let $1\le q\le\infty$ and $0\le\ell\le K$ with $\ell<\theta+N(1-1/q)$. 
Assume $p>1+\theta/N$ and 
put $F(u(x,t)):=\lambda|u(x,t)|^{p-1}u(x,t)$. 
Assertion~(a) follows from the similar argument to that of the proof of Corollary~\ref{Corollary:5.1}. 
Furthermore, for any $\sigma>0$, by \eqref{eq:2.2} and \eqref{eq:5.12} we have
\begin{equation*}
	\begin{aligned}
 	& t^{\frac{N}{\theta}\left(1-\frac{1}{q}\right)-\frac{\ell}{\theta}}|||U_0(t)-M_*g(t)|||_{q,\ell}\\
 	& \le C\int_t^\infty |M_0(F(s))|\,ds+C\sum_{1\le |\alpha|\le K}t^{-\frac{|\alpha|}{\theta}}|M_\alpha(\varphi)|
	+C\sum_{m=1}^{K_\theta}\sum_{|\alpha|\le K}t^{-m-\frac{|\alpha|}{\theta}}|M_\alpha(\varphi)|\\
	&\qquad\qquad
	+C\sum_{m=1}^{K_\theta}
	\sum_{|\alpha|\le K}t^{-m-\frac{|\alpha|}{\theta}}\int_0^t (s+1)^m |M_\alpha(F(s))|\,ds\\
 	& \qquad\qquad\qquad\qquad
 	+C\sum_{1\le|\alpha|\le K}t^{-\frac{|\alpha|}{\theta}}\int_0^t |M_\alpha(F(s))|\,ds\\
 	& =O\left(t^{-(A_p-1)}\right)+O(t^{-\frac{1}{\theta}})+O(t^{-1})+O\left(t^{-1}\int_1^t s^{1-A_p}\,ds\right)
	+O\left(t^{-\frac{1}{\theta}}\int_1^t s^{\frac{1}{\theta}-A_p}\,ds\right)\\
 	& =O\left(t^{-(A_p-1)+\sigma}\right)+O(t^{-1})+O(t^{-\frac{1}{\theta}})
	=O(h_\sigma(t))
	\end{aligned}
\end{equation*}
as $t\to\infty$. 
This together with \eqref{eq:5.5} implies that 
\begin{equation}
\label{eq:5.17}
	t^{\frac{N}{\theta}\left(1-\frac{1}{q}\right)-\frac{\ell}{\theta}}|||u(t)-M_*g(t)|||_{q,\ell}
	=o\left(t^{-\frac{K}{\theta}}\right)+O(h_\sigma(t))
	\quad\mbox{as}\quad t\to\infty.
\end{equation}
Furthermore, 
combining Lemma~\ref{Lemma:5.1} and \eqref{eq:5.17}, we have 
\begin{equation}
\label{eq:5.18}
	E_{K,q}[F(u)-F_\infty](t)=o\left(t^{-A_p-\frac{K}{\theta}}\right)+O\left(t^{-A_p}h_\sigma(t)\right)
\end{equation}
as $t\to\infty$.
On the other hand, it follows that
$$
	\begin{aligned}
	w(x,t) & :=u(x,t)-\int_0^t e^{-(t-s)(-\Delta)^{\frac{\theta}{2}}}F_\infty(s)\,ds\\
 	& \,\,=\left[e^{-t(-\Delta)^{\frac{\theta}{2}}}\varphi\right](x)
	+\int_0^t e^{-(t-s)(-\Delta)^{\frac{\theta}{2}}}[F(u(s))-F_\infty(s)]\,ds.
	\end{aligned}
$$
Applying Theorem~\ref{Theorem:1.1}, 
for any $\varepsilon>0$ and $T>0$, we obtain
\begin{equation}
\label{eq:5.19}
	\begin{aligned}
 	& t^{\frac{N}{\theta}\left(1-\frac{1}{q}\right)-\frac{\ell}{\theta}}|||w(t)-w_*(t)|||_{q,\ell}\\
 	& \le\varepsilon t^{-\frac{K}{\theta}}
	+C_*t^{-\frac{K}{\theta}}\int_T^t s^{\frac{K}{\theta}}E_{K,q}[F(u)-F_\infty](s)\,ds
	\end{aligned}
\end{equation}
as $t\to\infty$, where $C_*$ is a positive constant independent of $\varepsilon$ and $T$ and
$$
	w_*(x,t)=\sum_{m=0}^{K_\theta}\sum_{|\alpha|\le K}
	\left(M_\alpha(\varphi)+\int_0^t (s+1)^mM_\alpha(F(u(s))-F_\infty(s))\,ds\right) g_{\alpha,m}(x,t).
$$
Since $\varepsilon$ is arbitrary, 
by \eqref{eq:5.18} and \eqref{eq:5.19} we see that
$$
	t^{\frac{N}{\theta}\left(1-\frac{1}{q}\right)-\frac{\ell}{\theta}}|||w(t)-w_*(t)|||_{q,\ell}
	=o\left(t^{-\frac{K}{\theta}}\right)+O\left(t^{-\frac{K}{\theta}}
	\int_T^t s^{\frac{K}{\theta}-A_p}h_\sigma(s)\,ds\right)
$$
as $t\to\infty$. 
This implies assertion~(b). 
Thus Theorem~\ref{Theorem:1.2} follows. 
$\Box$
\begin{remark}
\label{Remark:5.2}
	Let $u$ be a solution to the Cauchy problem for a nonlinear fractional diffusion equation 
	and possess the mass conservation law, that is, 
	$\int_{{\mathbb R}^N}u(x,t)\,dx$ is independent of $t$. 
	The mass conservation law has often played an important role 
	in the study of HOAE of solutions to various nonlinear problems, 
	see e.g. \cites{DC, EZ, FM, NY, YM, YS01, YS02, YS03, YT1, YT2}. 
	Then the arguments in the proof of Theorem~{\rm\ref{Theorem:4.1}} are valid 
	for the Cauchy problem.
\end{remark}
{\bf Acknowledgements.} 
The authors of this paper were supported in part by JSPS KAKENHI Grant Number JP19H05599. 
The second author was also supported in part by JSPS KAKENHI Grant Number JP20K03689. 
  

\begin{bibdiv}
\begin{biblist}
\bib{AJY}{article}{
   author={Achleitner, Franz},
   author={J\"{u}ngel, Ansgar},
   author={Yamamoto, Masakazu},
   title={Large-time asymptotics of a fractional drift-diffusion-Poisson
   system via the entropy method},
   journal={Nonlinear Anal.},
   volume={179},
   date={2019},
   pages={270--293},
}
\bib{AF}{article}{
   author={Amann, H.},
   author={Fila, M.},
   title={A Fujita-type theorem for the Laplace equation with a dynamical
   boundary condition},
   journal={Acta Math. Univ. Comenian. (N.S.)},
   volume={66},
   date={1997},
   pages={321--328},
}
\bib{BT}{article}{
   author={Bogdan, Krzysztof},
   author={Jakubowski, Tomasz},
   title={Estimates of heat kernel of fractional Laplacian perturbed by
   gradient operators},
   journal={Comm. Math. Phys.},
   volume={271},
   date={2007},
   pages={179--198},
}
\bib{BK}{article}{
   author={Brandolese, Lorenzo},
   author={Karch, Grzegorz},
   title={Far field asymptotics of solutions to convection equation with
   anomalous diffusion},
   journal={J. Evol. Equ.},
   volume={8},
   date={2008},
   pages={307--326},
}
\bib{DC}{article}{
   author={Duro, Gema},
   author={Carpio, Ana},
   title={Asymptotic profiles for convection-diffusion equations with
   variable diffusion},
   journal={Nonlinear Anal.},
   volume={45},
   date={2001},
   pages={407--433},
}
\bib{EZ}{article}{
   author={Escobedo, Miguel},
   author={Zuazua, Enrike},
   title={Large time behavior for convection-diffusion equations in ${\mathbb R}^N$},
   journal={J. Funct. Anal.},
   volume={100},
   date={1991},
   pages={119--161},
}
\bib{FIK}{article}{
   author={Fila, Marek},
   author={Ishige, Kazuhiro},
   author={Kawakami, Tatsuki},
   title={Convergence to the Poisson kernel for the Laplace equation with a
   nonlinear dynamical boundary condition},
   journal={Commun. Pure Appl. Anal.},
   volume={11},
   date={2012},
   pages={1285--1301},
}
\bib{FK}{article}{
   author={Fino, Ahmad},
   author={Karch, Grzegorz},
   title={Decay of mass for nonlinear equation with fractional Laplacian},
   journal={Monatsh. Math.},
   volume={160},
   date={2010},
   pages={375--384},
}
\bib{FM}{article}{
   author={Fujigaki, Yoshiko},
   author={Miyakawa, Tetsuro},
   title={Asymptotic profiles of nonstationary incompressible Navier-Stokes
   flows in the whole space},
   journal={SIAM J. Math. Anal.},
   volume={33},
   date={2001},
   pages={523--544},
}
\bib{HKN}{article}{
   author={Hayashi, Nakao},
   author={Kaikina, Elena I.},
   author={Naumkin, Pavel I.},
   title={Asymptotics for fractional nonlinear heat equations},
   journal={J. London Math. Soc. (2)},
   volume={72},
   date={2005},
   pages={663--688},
}
\bib{IIK}{article}{
   author={Ishige, Kazuhiro},
   author={Ishiwata, Michinori},
   author={Kawakami, Tatsuki},
   title={The decay of the solutions for the heat equation with a potential},
   journal={Indiana Univ. Math. J.},
   volume={58},
   date={2009},
   pages={2673--2707},
}
\bib{IK}{article}{
   author={Ishige, Kazuhiro},
   author={Kawakami, Tatsuki},
   title={Refined asymptotic profiles for a semilinear heat equation},
   journal={Math. Ann.},
   volume={353},
   date={2012},
   pages={161--192},
   issn={0025-5831},
}
\bib{IKK01}{article}{
   author={Ishige, Kazuhiro},
   author={Kawakami, Tatsuki},
   author={Kobayashi, Kanako},
   title={Global solutions for a nonlinear integral equation with a
   generalized heat kernel},
   journal={Discrete Contin. Dyn. Syst. Ser. S},
   volume={7},
   date={2014},
   pages={767--783},
}
\bib{IKK02}{article}{
   author={Ishige, Kazuhiro},
   author={Kawakami, Tatsuki},
   author={Kobayashi, Kanako},
   title={Asymptotics for a nonlinear integral equation with a generalized
   heat kernel},
   journal={J. Evol. Equ.},
   volume={14},
   date={2014},
   pages={749--777},
}
\bib{IKM}{article}{
   author={Ishige, Kazuhiro},
   author={Kawakami, Tatsuki},
   author={Michihisa, Hironori},
   title={Asymptotic expansions of solutions of fractional diffusion
   equations},
   journal={SIAM J. Math. Anal.},
   volume={49},
   date={2017},
   pages={2167--2190},
}
\bib{IKO}{article}{
   author={Ishige, Kazuhiro},
   author={Kawakami, Tatsuki},
   author={Okabe, Shinya},
   title={Existence of solutions to nonlinear parabolic equations via majorant integral kernel}, 
   journal={preprint (arXiv:2101.06581)},
}
\bib{Iwa}{article}{
   author={Iwabuchi, Tsukasa},
   title={Global solutions for the critical Burgers equation in the Besov
   spaces and the large time behavior},
   journal={Ann. Inst. H. Poincar\'{e} Anal. Non Lin\'{e}aire},
   volume={32},
   date={2015},
   pages={687--713},
}
\bib{NY}{article}{
   author={Nagai, Toshitaka},
   author={Yamada, Tetsuya},
   title={Large time behavior of bounded solutions to a parabolic system of
   chemotaxis in the whole space},
   journal={J. Math. Anal. Appl.},
   volume={336},
   date={2007},
   pages={704--726},
}
\bib{S}{article}{
   author={Sugitani, Sadao},
   title={On nonexistence of global solutions for some nonlinear integral
   equations},
   journal={Osaka Math. J.},
   volume={12},
   date={1975},
   pages={45--51},
}
\bib{W}{article}{
   author={Weissler, Fred B.},
   title={Existence and nonexistence of global solutions for a semilinear
   heat equation},
   journal={Israel J. Math.},
   volume={38},
   date={1981},
   pages={29--40},
}
\bib{YM}{article}{
   author={Yamamoto, Masakazu},
   title={Asymptotic expansion of solutions to the dissipative equation with
   fractional Laplacian},
   journal={SIAM J. Math. Anal.},
   volume={44},
   date={2012},
   pages={3786--3805; {\it Erratum} in SIAM J. Math. Anal. {\bf 48} (2016), 3037--3038},
}
\bib{YS01}{article}{
   author={Yamamoto, Masakazu},
   author={Sugiyama, Yuusuke},
   title={Asymptotic expansion of solutions to the drift-diffusion equation
   with fractional dissipation},
   journal={Nonlinear Anal.},
   volume={141},
   date={2016},
   pages={57--87},
}
\bib{YS02}{article}{
   author={Yamamoto, Masakazu},
   author={Sugiyama, Yuusuke},
   title={Spatial-decay of solutions to the quasi-geostrophic equation with
   the critical and supercritical dissipation},
   journal={Nonlinearity},
   volume={32},
   date={2019},
   pages={2467--2480},
}
\bib{YS03}{article}{
   author={Yamamoto, Masakazu},
   author={Sugiyama, Yuusuke},
   title={Optimal estimates for far field asymptotics of solutions to the
   quasi-geostrophic equation},
   journal={Proc. Amer. Math. Soc.},
   volume={149},
   date={2021},
   pages={1099--1110},
}
\bib{YT1}{article}{
   author={Yamada, Tetsuya},
   title={Moment estimates and higher-order asymptotic expansions of
   solutions to a parabolic system of chemotaxis in the whole space},
   journal={Funkcial. Ekvac.},
   volume={54},
   date={2011},
   pages={15--51},
}
\bib{YT2}{article}{
   author={Yamada, T.},
   title={Improvement of convergence rates for a parabolic system of
   chemotaxis in the whole space},
   journal={Math. Methods Appl. Sci.},
   volume={34},
   date={2011},
   pages={2103--2124},
}
\end{biblist}
\end{bibdiv}
\end{document}